\pdfoutput=1
\RequirePackage{ifpdf}
\ifpdf 
\documentclass[pdftex]{sigma}
\else
\documentclass{sigma}
\fi

\numberwithin{equation}{section}

\newtheorem{thm}{Theorem}[section]
\newtheorem{coro}[thm]{Corollary}

\newtheorem{problem}[thm]{Problem}

\theoremstyle{definition}
\newtheorem{dfn}[thm]{Definition}

\begin{document}


\renewcommand{\thefootnote}{$\star$}

\newcommand{\arXivNumber}{1604.03133}

\renewcommand{\PaperNumber}{078}

\FirstPageHeading

\ShortArticleName{Universality Limits}
\ArticleName{An Update on Local Universality Limits\\ for Correlation Functions Generated\\ by Unitary Ensembles\footnote{This paper is a~contribution to the Special Issue on Asymptotics and Universality in Random Matrices, Random Growth Processes, Integrable Systems and Statistical Physics in honor of Percy Deift and Craig Tracy. The full collection is available at \href{http://www.emis.de/journals/SIGMA/Deift-Tracy.html}{http://www.emis.de/journals/SIGMA/Deift-Tracy.html}}}

\Author{Doron S.~LUBINSKY}
\AuthorNameForHeading{D.S.~Lubinsky}
\Address{School of Mathematics, Georgia Institute of Technology, Atlanta, GA 30332-0160 USA}
\Email{\href{mailto:lubinsky@math.gatech.edu}{lubinsky@math.gatech.edu}}

\ArticleDates{Received April 05, 2016, in f\/inal form August 05, 2016; Published online August 10, 2016}

\Abstract{We survey the current status of universality limits for $m$-point correlation functions in the bulk and at the edge for unitary ensembles, primarily when the limiting kernels are Airy, Bessel, or Sine kernels. In particular, we consider underlying measures on compact intervals, and f\/ixed and varying
exponential weights, as well as universality limits for a~variety of orthogonal systems. The scope of the survey is quite narrow: we do not consider~$\beta$ ensembles for $\beta \neq 2$, nor general Hermitian matrices with independent entries, let alone more general settings. We include some open problems.}

\Classification{15B52; 60B20; 60F99; 42C05; 33C50}

\renewcommand{\thefootnote}{\arabic{footnote}}
\setcounter{footnote}{0}

\section{Introduction}

The remarkable connections between random matrices and other topics are clearly exposited in many texts, so we shall not discuss it in this brief review. Readers may consult \cite{Akemanetal2011,Andersonetal2010,Baiketal2008,BleherIts2001,Deift1999,DeiftGioev2009,Forrester2010,Mehta2004,Tao2012}. The 2003 short
historical survey of Forrester, Snaith and Verbaarschot \cite{Forresteretal2003} is still very useful. We simply start with a probability distribution $P^{(n) }$ on the space of $n$ by~$n$ Hermitian matrices $M=( t_{ij}) _{1\leq i,j\leq n}$:
\begin{gather*}
P^{(n) }(M) =\frac{1}{Z_{n}}w(M) dM =\frac{1}{Z_{n}}w(M) \left(\prod_{j=1}^{n}dt_{jj}\right) \left( \prod_{j<k}d\left(\operatorname{Re}t_{jk}\right) d\left( \operatorname{Im}t_{jk}\right) \right).
\end{gather*}
Here $w$ is some non-negative function def\/ined on Hermitian matrices, and $Z_{n}$ is a normalizing constant, there to ensure that $P^{(n)}$ is indeed a~probability distribution. The most important case is
\begin{gather*}
w(M) =\exp ( -2n\operatorname{tr}Q(M)),
\end{gather*}
for appropriate functions $Q$. In particular, the choice $Q(M)
=M^{2}$, leads to the Gaussian unitary ensemble (apart from scaling) that
was considered by Wigner. When expressed in spectral form, that is, as a
probability distribution on the eigenvalues $x_{1}\leq x_{2}\leq \cdots\leq
x_{n}$ of $M$, it takes the form
\begin{gather*}
P^{(n) }( x_{1},x_{2},\dots,x_{n}) dx_{1}\cdots dx_{n}=\frac{1}{Z_{n}}\left( \prod\limits_{1\leq i<j<n}( x_{i}-x_{j})
^{2}\right) \left( \prod\limits_{j=1}^{n}w (x_{j}) dx_{j}\right).
\end{gather*}
See \cite[p.~102~f\/f.]{Deift1999}. Again, $Z_{n}$ is a normalizing constant, often called the partition function (constant). Note that $w$ now can be any
non-negative measurable function. We could also replace the absolutely continuous measure~$w(x) dx$ by $d\mu (x)$, where $\mu$ is a positive measure on the real line, but for the moment focus on absolutely continuous measures.

In most applications, we want to let $n\rightarrow \infty $, and obviously the $n$-fold density complicates issues. So we often integrate out most variables, forming marginal distributions. One particularly important quantity is the $m$-point correlation function \cite[p.~112]{Deift1999}, which in the following form includes a factor of $w(x_{1})w(x_{2})\cdots w(x_{m})$:
\begin{gather*}
R_{m,n}( w;x_{1},x_{2},\dots,x_{m})\nonumber \\
\qquad{}=\frac{n!}{(n-m)!}\frac{1}{Z_{n}}\int \cdots\int \left(\prod\limits_{1\leq i<j<n}( x_{i}-x_{j}) ^{2}\right) \left(\prod\limits_{j=1}^{n}w(x_{j}) \right) dx_{m+1}dx_{m+2}\cdots dx_{n}.
\end{gather*}
Typically, we f\/ix $m$, and study $R_{m,n}$ as $n\rightarrow \infty $. $R_{m,n}$ is useful in examining spacing of eigen\-values, and counting the expected number of eigenvalues in some set. For example, if $B$ is a~measu\-rable subset of~$\mathbb{R}$,
\begin{gather*}
\int_{B}\cdots \int_{B}R_{m,n}(w;x_{1},x_{2},\dots,x_{m}) dx_{1}dx_{2}\cdots dx_{m}
\end{gather*}
counts the expected number of $m$-tuples $(x_{1},x_{2},\dots,x_{m})$ of distinct eigenvalues with each $x_{j}\in B$.

There are several types of universality limits involving the $m$-point correlation function $R_{m,n}$. We shall examine those involving local or ``microlocal'' scaling about a point inside the interior of the support of~$w$ (the ``bulk'') or at the boundary of the support of $w$ (the ``edge''). When we deal with a sequence of weights, the so-called ``varying weights'' case, the edge might be a~moving, or ``soft edge'', or possibly a f\/ixed, or ``hard edge''.

Let us illustrate these through the examples of the Jacobi and Hermite weights. The Jacobi weight is def\/ined by
\begin{gather*}
w_{\alpha,\beta }(x) =(1-x)^{\alpha }(1+x) ^{\beta }, \qquad x\in (-1,1),
\end{gather*}
where $\alpha,\beta >-1$. Given $x\in (-1,1) $, $m\geq 1$, and real numbers $u_{1},u_{2},\dots,u_{m}$, its universality limit in the bulk takes the form (cf.\ \cite[Theorem~1.1, p.~1577]{KuijlaarsVanlessen2002})
\begin{gather*}
\lim_{n\rightarrow \infty }\!\left(\! \frac{\pi \sqrt{1-x^{2}}}{n}\!\right)^{m}\!\! R_{m,n}\!\left(\!w_{\alpha,\beta };x+\frac{u_{1}\pi \sqrt{1-x^{2}}}{n},x+\frac{u_{2}\pi \sqrt{1-x^{2}}}{n},\dots,
 x+\frac{u_{m}\pi \sqrt{1-x^{2}}}{n}\right) \\
\qquad {}=\det (\mathbb{S}( u_{i}-u_{j})) _{1\leq i,j\leq m}.
\end{gather*}
Here
\begin{gather*}
\mathbb{S}(x) =\frac{\sin \pi x}{\pi x}
\end{gather*}
is the sinc kernel, and the factor $\pi \sqrt{1-x^{2}}$ is the reciprocal of the equilibrium density for $(-1,1)$, see Section~\ref{section2}. At the
``hard edge'' $x=1$, we instead have for $u_{1},u_{2},\dots,u_{m}>0$ (cf.\ \cite[Theorem~1.1, p.~1577]{KuijlaarsVanlessen2002})
\begin{gather*}
\lim_{n\rightarrow \infty }\left( \frac{1}{2n^{2}}\right) ^{m}R_{m,n}\left(w_{\alpha,\beta };1-\frac{u_{1}}{2n^{2}},1-\frac{u_{2}}{2n^{2}},\dots,1-\frac{u_{m}}{2n^{2}}\right) =\det \left( \mathbb{J}_{\alpha }\left(u_{i},u_{j}\right) \right)_{1\leq i,j\leq m}.
\end{gather*}
Here, $\mathbb{J}_{\alpha }$ denotes the Bessel kernel of order $\alpha$, involving the usual Bessel function $J_{\alpha }$ of order $\alpha$,
\begin{gather*}
\mathbb{J}_{\alpha }(u,v) =\frac{J_{\alpha }\left( \sqrt{u}\right) \sqrt{v}J_{\alpha }^{\prime }\left( \sqrt{v}\right) -J_{\alpha
}\left( \sqrt{v}\right) \sqrt{u}J_{\alpha }^{\prime }\left( \sqrt{u}\right) }{2(u-v)}.
\end{gather*}

For the varying Hermite weight
\begin{gather*}
w_{H}^{2n}(x) =\exp \big({-}2nx^{2}\big), \qquad x\in (-\infty,\infty),
\end{gather*}
the ``bulk'' is not the whole real line, but the f\/inite interval $[-1,1] $~-- it is here where the relevant equilibrium density is supported (more on this in Section 5). The universality limit in the bulk takes the form (cf.\ \cite[p.~253, Theorem~8.16]{Deift1999}, \cite[p.~1348, Theorem~1.4]{Deiftetal1999}, \cite[p.~283]{Forrester2010})
\begin{gather*}
\lim_{n\rightarrow \infty }\left( \frac{\pi }{2n\sqrt{1-x^{2}}}\right)
^{m}R_{m,n}\left( w_{H}^{2n};x+\frac{u_{1}\pi }{2n\sqrt{1-x^{2}}},x+\frac{u_{2}\pi }{2n\sqrt{1-x^{2}}},\dots,x+\frac{u_{m}\pi }{2n\sqrt{1-x^{2}}}\right)
\\
\qquad{} =\det \left( \mathbb{S}\left( u_{i}-u_{j}\right) \right) _{1\leq i,j\leq m},
\end{gather*}
and holds uniformly for $x$ in compact subsets of $(-1,1) $ and $u_{1},u_{2},\dots,u_{m}$ in compact subsets of the real line. The function $\frac{2}{\pi }\sqrt{1-x^{2}}$ is the equilibrium density for the external f\/ield $x^{2}$, see Section~\ref{section5}. The universality at the soft edge~$1$ takes
the form \cite[equation~(1.6), p.~4]{LevinLubinsky2010}, \cite[p.~152]{TracyWidom1994}
\begin{gather*}
\lim_{n\rightarrow \infty }\! \left( \frac{1}{2n^{2/3}}\right)
^{m}\!\! R_{m,n}\!\left( w_{H}^{2n};1+\frac{u_{1}}{2n^{2/3}},1+\frac{u_{2}}{2n^{2/3}
},\dots,1+\frac{u_{m}}{2n^{2/3}}\right) =\det \left( \mathbb{A}{\rm i}( u_{i},u_{j}) \right) _{1\leq i,j\leq m},
\end{gather*}
where $\mathbb{A}{\rm i}(\cdot,\cdot) $ is the Airy kernel, def\/ined by
\begin{gather*}
\mathbb{A}{\rm i}(u,v) = \begin{cases}
\dfrac{\operatorname{Ai}(u) \operatorname{Ai}^{\prime }(v) -\operatorname{Ai}^{\prime }(u) \operatorname{Ai}(v) }{u-v}, & u\neq v, \vspace{1mm}\\
\operatorname{Ai}^{\prime }(u) ^{2}-u\operatorname{Ai}(u) ^{2}, & u=v,
\end{cases}
\end{gather*}
and $\operatorname{Ai}$ is the Airy function, def\/ined on the real line by \cite[p.~53]{Olver1974}
\begin{gather*}
\operatorname{Ai}(x) =\frac{1}{\pi }\int_{0}^{\infty }\cos \left( \frac{1}{3}t^{3}+xt\right) dt.
\end{gather*}

The edge above is called ``soft'' because when we consider the non-scaled Hermite weight
\begin{gather*}
w_{H}(x) =e^{-x^{2}}, \qquad x\in (-\infty,\infty),
\end{gather*}
then the bulk becomes the growing sequence of intervals $\left( -\sqrt{2n},\sqrt{2n}\right) $ and $\sqrt{2n}$ is soft in the sense that it shifts~-- this terminology carries over to the point~$1$ after scaling. Now let us consider this f\/ixed Hermite weight~$w_{H}$. Let $\varepsilon \in (0,1)$. The universality limit in the bulk takes the form \cite[p.~4]{LevinLubinsky2010}
\begin{gather*}
\lim_{n\rightarrow \infty }\left( \frac{\pi }{\sqrt{2n-x^{2}}}\right)
^{m}R_{m,n}\left( w_{H};x+\frac{u_{1}\pi }{\sqrt{2n-x^{2}}},x+\frac{u_{2}\pi
}{\sqrt{2n-x^{2}}},\dots,x+\frac{u_{m}\pi }{\sqrt{2n-x^{2}}}\right) \\
\qquad{}=\det \left( \mathbb{S}\left( u_{i}-u_{j}\right) \right) _{1\leq i,j\leq m},
\end{gather*}
uniformly for $x$ in $\left( -\sqrt{2n}\left( 1-\varepsilon \right),\sqrt{2n} ( 1-\varepsilon) \right) $, and $u_{1},u_{2},\dots,u_{m}$ in compact subsets of the real line. At the soft-edge, the limit takes the form \cite[p.~286]{Forrester2010}, \cite[p.~4]{LevinLubinsky2010}
\begin{gather*}
\lim_{n\rightarrow \infty }\left( \frac{1}{\sqrt{2}n^{1/6}}\right)
^{m}R_{m,n}\left( w_{H};\sqrt{2n}\left( 1+\frac{u_{1}}{2n^{2/3}}\right),
\sqrt{2n}\left( 1+\frac{u_{2}}{2n^{2/3}}\right), \dots,\sqrt{2n}\left( 1+\frac{u_{m}}{2n^{2/3}}\right) \right) \\
\qquad{} =\det ( \mathbb{A}{\rm i} ( u_{i},u_{j} ) ) _{1\leq i,j\leq m}.
\end{gather*}
Of course universality limits for Laguerre weights $x^{\alpha }e^{-x}$ on $[0,\infty )$, $\alpha >-1,$ are closely related to those of Hermite weights.
There is a new feature, however~-- $0$ is a hard edge, the bulk becomes $\left( \varepsilon,4n( 1-\varepsilon) \right) $, and the soft edge is around~$4n$.

Do these perhaps overly technical limits deter the reader? We gathered them in one place for the reason that they are dispersed in the literature. The good news is that at least the limits in the bulk can be put into a unifying form, and we shall do that shortly.

Who f\/irst discovered the identity, involving orthonormal polynomials, that played a crucial role in establishing these limits? For a long time, I was under the impression that it was Gaudin and Mehta, but I was informed that I~was mistaken at Percy's 70th birthday conference, and that it was Freeman Dyson. If you read the 2003 survey of Forrester, Snaith and Verbaarschot, they credit the orthogonal polynomial method to Mehta, citing a 1960 technical report~\cite{Mehta1960B}. McLaughlin and Miller~\cite{McLaughlinMiller2006} refer to a 1960 paper of Mehta and Gaudin~\cite{MehtaGaudin1960}. I am grateful to Thomas Bothner~\cite{Bothner2016} for the following insight: the orthogonal polynomial method was almost certainly born in those papers of Mehta and Gaudin from 1960 \cite{Mehta1960,MehtaGaudin1960}. However, Freeman Dyson appears to be the f\/irst to use it in the computation of the $m$-point correlation function, in a 1962 paper~\cite{Dyson1962}.

Given a positive measure $\mu $ on the real line, with all f\/inite power moments $\int x^{j}d\mu (x) $, $j=0,1,2,\dots $,
\begin{gather*}
p_{n} ( \mu,x ) =\gamma _{n}x^{n}+\cdots, \qquad \gamma _{n}>0,
\end{gather*}
denotes its $n$th orthonormal polynomial, so that
\begin{gather*}
\int p_{n}(\mu,x) p_{m}(\mu,x) d\mu(x) =\delta _{mn}.
\end{gather*}
They satisfy the three term recurrence relation
\begin{gather}\label{eq1.6}
xp_{n}(\mu,x) =a_{n+1}p_{n+1}(\mu,x) +b_{n}p_{n}(\mu,x) +a_{n}p_{n-1}(\mu,x),\qquad n\geq 0,
\end{gather}
where $a_{n}=\frac{\gamma _{n-1}}{\gamma _{n}}>0$ and $b_{n}$ is real. Conversely any sequence of polynomials satisfying such a~three term recurrence is a~sequence of orthogonal polynomials for some positive measure on the real line. Def\/ine the $n$th reproducing kernel
\begin{gather*}
K_{n}(\mu,x,y) =\sum_{j=0}^{n-1}p_{j}(\mu,x) p_{j}(\mu,y)
\end{gather*}
and its normalized cousin
\begin{gather*}
\tilde{K}_{n}(\mu,x,y) =\mu ^{\prime }(x) ^{1/2}\mu ^{\prime }(y) ^{1/2}K_{n}(\mu,x,y).
\end{gather*}
Here $\mu ^{\prime }$ denotes the Radon--Nikodym derivative of the absolutely continuous component of $\mu $, that exists a.e.\ It is taken as~$0$ elsewhere. When $\mu $ is absolutely continuous and
\begin{gather*}
w=\mu ^{\prime },
\end{gather*}
we use the notation $p_{n}(w,x)$, $K_{n}(w,x,y) $, $\tilde{K}_{n}(w,x,y) $. The remarkable identity referred to above, asserts that the $m$-point correlation function for the weight~$w$ is given by \cite[p.~112]{Deift1999}
\begin{gather*}
R_{m,n} ( w;x_{1},x_{2},\dots,x_{m} ) =\det ( \tilde{K}_{n} ( w,x_{i},x_{j} )) _{1\leq i,j\leq m}.
\end{gather*}

Since the determinant has f\/ixed size, and the entries all have similar form, we see that the universality limits above reduce to the following asymptotics for the reproducing kernels. We omit the restrictions on the parameters, they are similar to those above:
\begin{enumerate}\itemsep=0pt
\item[(I)] \emph{The Jacobi weight $w_{\alpha,\beta }(x) =(1-x) ^{\alpha }(1+x) ^{\beta }$ on $(-1,1) $.}
\begin{enumerate}\itemsep=0pt
\item[(a)] \emph{In the bulk $(-1,1) $}
\begin{gather*}
\lim_{n\rightarrow \infty }\frac{\pi \sqrt{1-x^{2}}}{n}\tilde{K}_{n}\left(
w_{\alpha,\beta },x+\frac{\pi \sqrt{1-x^{2}}}{n}u,x+\frac{\pi \sqrt{1-x^{2}}}{n}v\right) =\mathbb{S}(u-v).
\end{gather*}
\item[(b)] \emph{At the hard edge $1$}
\begin{gather*}
\lim_{n\rightarrow \infty }\frac{1}{2n^{2}}\tilde{K}_{n}\left( w_{\alpha,\beta },1-\frac{u}{2n^{2}},1-\frac{v}{2n^{2}}\right) =\mathbb{J}_{\alpha
}(u,v).
\end{gather*}
\end{enumerate}

\item[(II)] \emph{The fixed Hermite weight $w_{H}(x) =e^{-x^{2}}$ on $( -\infty,\infty ) $.}
\begin{enumerate}\itemsep=0pt
\item[(a)] \emph{In the bulk $\left( -\sqrt{2n},\sqrt{2n}\right)$}
\begin{gather*}
\lim_{n\rightarrow \infty }\frac{\pi }{\sqrt{2n-x^{2}}}\tilde{K}_{n}\left(w_{H},x+\frac{\pi }{\sqrt{2n-x^{2}}}u,x+\frac{\pi }{\sqrt{2n-x^{2}}}v\right)
=\mathbb{S}(u-v).
\end{gather*}
\item[(b)] \emph{At the soft edge $\sqrt{2n}$}
\begin{gather*}
\lim_{n\rightarrow \infty }\frac{1}{\sqrt{2}n^{1/6}}\tilde{K}_{n}\left(w_{H},\sqrt{2n}\left( 1+\frac{u}{2n^{2/3}}\right),\sqrt{2n}\left( 1+\frac{v}{2n^{2/3}}\right) \right) =\mathbb{A}{\rm i}(u,v).
\end{gather*}
\end{enumerate}

\item[(III)] \emph{The varying Hermite weights $w_{H}^{2n}(x)=e^{-2nx^{2}}$ on $( -\infty,\infty) $.}
\begin{enumerate}\itemsep=0pt
\item[(a)] \emph{In the bulk $(-1,1)$}
\begin{gather*}
\lim_{n\rightarrow \infty }\frac{\pi }{2n\sqrt{1-x^{2}}}\tilde{K}_{n}\left(
w_{H}^{2n},x+\frac{\pi }{2n\sqrt{1-x^{2}}}u,x+\frac{\pi }{2n\sqrt{1-x^{2}}}v\right) =\mathbb{S}(u-v).
\end{gather*}
\item[(b)] \emph{\emph{At the soft edge $1$}}
\begin{gather*}
\lim_{n\rightarrow \infty }\frac{1}{2n^{2/3}}\tilde{K}_{n}\left(w_{H}^{2n},1+\frac{u}{2n^{2/3}},1+\frac{v}{2n^{2/3}}\right) =\mathbb{A}{\rm i}(u,v).
\end{gather*}
\end{enumerate}
\end{enumerate}

Fortunately, all the bulk limits in (I), (II), (III) may be recast in a~unif\/ied form, which we deliberately formulate vaguely:

\begin{enumerate}\itemsep=0pt
\item[(IV)] \emph{General bulk universality.} Let $\left\{ \mu _{n}\right\} $ be a sequence of measures, and $x$ lie in the bulk. Then
\begin{gather}\label{eq1.17}
\lim_{n\rightarrow \infty }\frac{\tilde{K}_{n}\left( \mu _{n},x+\frac{u}{\tilde{K}_{n}(\mu _{n},x,x) },x+\frac{v}{\tilde{K}_{n}(\mu_{n},x,x) }\right) }{\tilde{K}_{n}(\mu _{n},x,x) }=\mathbb{S}(u-v).
\end{gather}
Similarly, one can do this for soft and hard edges, though the general formulation is less useful than in the bulk case, since soft and hard edge universality is inherently more special, and not so widely established.

\item[(V)] \emph{General soft edge universality.} If $1$ is a soft edge for $\mu _{n}$, $n\geq 1$, then \cite[p.~5]{LevinLubinsky2010}
\begin{gather}\label{eq1.18}
\lim_{n\rightarrow \infty }\frac{\tilde{K}_{n}\left( \mu _{n},1+\frac{\mathbb{A}{\rm i}(0,0) }{\tilde{K}_{n}(\mu_{n},1,1) }u,x+
\frac{\mathbb{A}{\rm i}(0,0) }{\tilde{K}_{n}(\mu_{n},1,1)
}v\right) }{\tilde{K}_{n}(\mu_{n},1,1) }=\frac{\mathbb{A}{\rm i}(u,v) }{\mathbb{A}{\rm i}(0,0) }.
\end{gather}
\item[(VI)] \emph{General hard edge universality.} If $1$ is a hard edge for $\mu _{n}$, $n\geq 1$, then \cite[p.~5]{Lubinsky2008D} for the Bessel kernel of order $\alpha $,
\begin{gather}
\lim_{n\rightarrow \infty }\frac{K_{n}\left( \mu _{n},1-\left( \frac{\mathbb{J}_{\alpha }^{\ast }(0,0) }{K_{n}(\mu_{n},1,1) }\right) ^{1/(\alpha +1) }u^{2},1-\left( \frac{\mathbb{J}_{\alpha }^{\ast }(0,0) }{K_{n}(\mu_{n},1,1) }\right)^{1/(\alpha +1) }v^{2}\right) }{K_{n}(\mu_{n},1,1) }\nonumber\\
\qquad{} =\frac{\mathbb{J}_{\alpha }^{\ast }\left( u^{2},v^{2}\right) }{\mathbb{J}_{\alpha }^{\ast }(0,0) },\label{eq1.19}
\end{gather}
where
\begin{gather*}
\mathbb{J}_{\alpha }^{\ast }(u,v) =\mathbb{J}_{\alpha } (u.v) /\big\{ u^{\alpha /2}v^{\alpha /2}\big\}.
\end{gather*}
The advantage of $\mathbb{J}_{\alpha }^{\ast }(u,v) $ over $\mathbb{J}_{\alpha }( u.v) $ is that the former is an entire function in two variables.
\end{enumerate}

In subsequent sections, we discuss
\begin{enumerate}\itemsep=0pt
\item[] \textit{How universal is universality?}
\end{enumerate}
That is, how general can the measures $\{ \mu _{n}\} $ in (\ref{eq1.17})--(\ref{eq1.19}), be?

We emphasize that we focus on a very narrow slice of universality limits for correlation functions. We omit the case of general $\beta $ $( \beta\neq 2) $ ensembles, of matrices with independent entries, double scaling limits, biorthogonal ensembles, \dots. Nor do we cover much the case of kernels other than the Airy, Bessel, or Sine Kernels, such as arise when equilibrium densities have zeros or inf\/inities inside the support of the equilibrium measure. We omit other universal features of eigenvalues, such as Gaussian behavior of local f\/luctuations of eigenvalues~\cite{Zhang2015} or of linear statistics~\cite{BreuerDuits2016}, and mesoscopic f\/luctuations~\cite{BreuerDuits2015}. Recent overviews of the more general case are given by L.~Erd\H{o}s in \cite{Erdos2011}, A.~Kuijlaars in \cite{Kuijlaars2011}, and T.~Tao and V.~Vu in \cite{TaoVu2014}.

The paper is structured as follows: in Section~\ref{section2}, we consider compactly supported measures. In Section~\ref{section3}, we consider applications of universality to orthogonal polynomial quantities. In Section~\ref{section4}, we consider universality for other orthogonal systems arising from a f\/ixed measure. In Section~\ref{section5}, we consider varying exponential weights, and in Section~\ref{section6}, f\/ixed exponential weights.

\section{Measures with compact support: universality}\label{section2}

In this section, we consider a (f\/ixed) measure $\mu $ with compact support, and f\/irst examine limits in the bulk. Much of my own research has dealt with this case, though referees have sometimes commented that this is not a case of physical interest. Nevertheless, as we shall see, it has intrinsic interest, and applications in orthogonal and random polynomials, especially in questions regarding zero distribution.

The most basic tool used is the Christof\/fel--Darboux formula,
\begin{gather*}
K_{n}(\mu,x,t) =\frac{\gamma _{n-1}}{\gamma _{n}}\frac{p_{n}(\mu,x) p_{n-1}(\mu,t) -p_{n-1}( \mu,x) p_{n}(\mu,t) }{x-t}.
\end{gather*}
If we have asymptotics for $p_{n}(\mu,x) $ as $n\rightarrow \infty $, such as
\begin{gather}\label{eq2.1}
\mu ^{\prime }(x) ^{1/2}p_{n}(\mu,x) \big(1-x^{2}\big) ^{1/4}=\sqrt{\frac{2}{\pi }}\cos ( n\theta +g(\theta )) +o(1),
\end{gather}
where $g$ is some continuous function, and $x=\cos \theta $, which is ``often'' true when $\mu $ is supported on $[-1,1] $, then substituting this into the Christof\/fel--Darboux formula with appropriate choices of~$x$,~$t$ yields a bulk limit, moduli some other minor factors.

One way to establish much more powerful asymptotics than \eqref{eq2.1} is the Riemann--Hilbert steepest descent method, pioneered by Deift and Zhou in the
1990's \cite{DeiftZhou1993, DeiftZhou1995}. It has revolutionized our understanding of asymptotics of orthogonal polynomials, as well as many other types of asymptotics. For generalized Jacobi type weights, it has been used by Arno Kuijlaars and his collaborators and students in a series of papers such as \cite{Kuijlaarsetal2004,KuijlaarsVanlessen2002}. Here \cite[Theorem~1.1, p.~1577]{KuijlaarsVanlessen2002} is one of their results:

\begin{thm} Let $h\colon [-1,1] \rightarrow ( 0,\infty ) $ be the restriction to $[-1,1] $ of a function analytic in a~neighborhood of $[-1,1] $. Let
$\alpha,\beta >-1$ and
\begin{gather*}
w(x) =h(x) (1-x) ^{\alpha }(1+x) ^{\beta }, \qquad x\in (-1,1).
\end{gather*}
Then uniformly for $x$ in compact subsets of $(-1,1) $ and $u$, $v$ in compact subsets of $\mathbb{R}$, we have
\begin{gather*}
\frac{\pi \sqrt{1-x^{2}}}{n}\tilde{K}_{n}\left( w,x+\frac{u\pi \sqrt{1-x^{2}}}{n},x+\frac{v\pi \sqrt{1-x^{2}}}{n}\right) =\mathbb{S}(u-v)
+O\left( \frac{1}{n}\right).
\end{gather*}
Moreover, uniformly for $u$, $v$ in compact subsets of $( 0,\infty) $,
\begin{gather*}
\frac{1}{2n^{2}}\tilde{K}_{n}\left( w,1-\frac{u}{2n^{2}},1-\frac{v}{2n^{2}}\right) =\mathbb{J}_{\alpha }(u,v) +O\left( \frac{u^{\alpha
/2}v^{\alpha /2}}{n}\right).
\end{gather*}
\end{thm}

As far as I know, these are still the best convergence rates for universality limits for Jacobi weights. The Riemann--Hilbert method actually yields these directly without having to substitute asymptotics into the Christof\/fel--Darboux formula.

Since pointwise asymptotics for $p_{n}(\mu,x) $ of the form~\eqref{eq2.1} typically require smoothness restrictions on~$\mu $, one needs new
ideas to establish universality for more general weights or measures. Inspired by Percy's 60th birthday conference, the author came up with a~comparison method: if $\mu $ and $\nu $ are positive measures with $\mu \leq \nu $, then for all real $x$, $y$, \cite[p.~919]{Lubinsky2009}
\begin{gather}\label{eq2.2}
\big\vert K_{n}(\mu,x,y) -K_{n}( \nu,x,y) \big\vert /K_{n}( \mu,x,x) \leq \left( \frac{K_{n}( \mu
,y,y) }{K_{n}( \mu,x,x) }\right) ^{1/2}\left[ 1-\frac{K_{n}( \nu,x,x) }{K_{n}( \mu,x,x) }\right] ^{1/2}.
\end{gather}
In particular, if $x$ and $y$ vary with $n$, and as $n\rightarrow \infty $, $\frac{K_{n}( \nu,x,x) }{K_{n}( \mu,x,x) }$ has limit~$1$, while $\frac{K_{n}(\mu,y,y) }{K_{n}(\mu,x,x) }$ remains bounded, then $K_{n}(\mu,x,y) $ and $K_{n}(\nu,x,y) $ have the same asymptotic. This inequality enables us to use universality limits for a larger ``nice'' measure $\nu $ to obtain the same for a ``not so nice'' measure~$\mu$, which is locally the same as~$\nu$.

We also need the concept of regularity of a measure in the sense of Stahl, Totik, and Ullmann (not to be confused with Borel regular measures). We say~$\mu$ is regular on $[-1,1] $ if
\begin{gather*}
\lim_{n\rightarrow \infty }\left( \sup_{\deg (P) \leq n}\frac{\Vert P\Vert _{L_{\infty }[-1,1] }}{\left( \int
\vert P\vert ^{2}d\mu \right) ^{1/2}}\right) ^{1/n}=1.
\end{gather*}
Thus sup norms of sequences of polynomials have the same $n$th root behavior as their $L_{2}(\mu) $ norms. This is known to be true if $\mu
^{\prime }>0$ a.e.\ in $[-1,1] $, though much less is required~\cite{StahlTotik1992}. An equivalent formulation involves the leading coef\/f\/icients $\{ \gamma _{n}\} $ of the orthonormal polynomials for~$\mu$:
\begin{gather*}
\lim_{n\rightarrow \infty }\gamma _{n}^{1/n}=\frac{1}{2}.
\end{gather*}

Using \eqref{eq2.2}, we proved \cite[Theorem~1.1, pp.~916--917]{Lubinsky2009}:

\begin{thm} \label{thm2.2} Let $\mu $ have support $[-1,1] $ and be regular. Let $x\in (-1,1) $ and assume $\mu $ is absolutely continuous in an open set containing~$x$. Assume moreover, that $\mu ^{\prime }$ is positive and continuous~at $x$. Then uniformly for $u$, $v$ in compact subsets of the real line, we have
\begin{gather*}
\lim_{n\rightarrow \infty }\frac{\widetilde{K}_{n}\left( \mu,x+\frac{u\pi
\sqrt{1-x^{2}}}{n},x+\frac{v\pi \sqrt{1-x^{2}}}{n}\right) }{\widetilde{K}_{n}(\mu,x,x) }=\mathbb{S}(u-v).
\end{gather*}
If the hypotheses hold in a compact set $J$, then the conclusion holds uniformly for $x\in J$.
\end{thm}

At the hard edge, the comparison technique yielded \cite[p.~283]{Lubinsky2008B}

\begin{thm}\label{thm2.3}Let $\mu $ be a finite positive Borel measure on $(-1,1) $ that is regular. Assume that for some $\rho >0$, $\mu $ is absolutely continuous in $J=[ 1-\rho,1]$, and in~$J$, its absolutely continuous component has the form $w(x) =h(x) (1-x) ^{\alpha }(1+x) ^{\beta }$, where $\alpha,\beta >-1$ and $h$ is a measurable function defined in $[1-\rho,1]$. Assume that $h(1) >0$ and $h$ is continuous at~$1$. Then uniformly for $u$, $v$ in compact subsets of $(0,\infty) $, we have
\begin{gather*}
\lim_{n\rightarrow \infty }\frac{1}{2n^{2}}\tilde{K}_{n}\left( 1-\frac{u}{2n^{2}},1-\frac{v}{2n^{2}}\right) =\mathbb{J}_{\alpha }(u,v).
\end{gather*}
If $\alpha \geq 0$, we may allow $u$, $v$ to lie in compact subsets of $[0,\infty )$.
\end{thm}

The real potential of the inequality \eqref{eq2.2} was soon explored by Findley, Simon and Totik \cite{Findley2008,Simon2007, Totik2009}. It was Findley who replaced continuity of $\mu ^{\prime }$ by the Szeg\H{o} condition on $[-1,1] $. Totik used the method of ``polynomial pullbacks'', which is based on the observation that if $P$ is a polynomial, then $P^{[-1] }[-1,1] $ consists of f\/initely many intervals. This allows one to pass from asymptotics for $[-1,1] $ to f\/initely many intervals. In turn, one can use the latter to approximate arbitrary compact sets. Barry Simon used instead Jost functions and dealt with the case of measures supported on several intervals.

To state Totik's result, we need a little more potential theory. Let $J$ be a compact subset of the real line, which we assume (for simplicity), has non-empty interior. We minimize the energy
\begin{gather*}
I[\nu] =\iint \log \frac{1}{\vert x-y\vert }d\nu (x) d\nu (y),
\end{gather*}
over all probability measures $\nu $ on $J$. The logarithmic capacity of $J$ is
\begin{gather*}
\operatorname{cap}(J) =\exp \left( -\inf_{\operatorname{supp}[\nu] \subset J,\nu (J) =1}I[\nu] \right).
\end{gather*}
There is a unique minimizing measure, called the equilibrium measure of~$J$. It need not be absolutely continuous, over all of $J$, but will be in the interior of~$J$ \cite[p.~216]{SaffTotik1997}. In the sequel, we represent it as $\omega (x) dx$ in the interior of~$J$, and call $\omega $ \textit{the equilibrium density of}~$J$. In the special case that $J=[a,b] $, then its capacity is $\frac{b-a}{4}$, and its equilibrium density is $\frac{1}{\pi \sqrt{(x-a)(b-x)}}$. In particular, for~$[-1,1] $, the equilibrium density is $\frac{1}{\pi \sqrt{1-x^{2}}}$, which is the underlying reason for this factor in the universality limits above.

We have already def\/ined the notion of regularity of a measure supported on $[-1,1] $. For $\mu $ with general compact support~$J$, regularity means that
\begin{gather*}
\lim_{n\rightarrow \infty }\gamma _{n}^{1/n}=\frac{1}{\operatorname{cap}(J) }.
\end{gather*}
An equivalent formulation is \cite[p.~66]{StahlTotik1992}
\begin{gather*}
\limsup_{n\rightarrow \infty }\left( \sup_{\deg (P) \leq n}\frac{\vert P(z)\vert ^{2}}{\int \vert P\vert ^{2}d\mu }\right) ^{1/n}\leq 1
\end{gather*}
for q.e.\ $z$ in $\operatorname{supp} [ \mu ] $, that is except possibly on a set of capacity~$0$. In particular, if~$J$ consists of f\/initely many intervals, it suf\/f\/ices that $\mu ^{\prime }>0$ a.e.\ on~$J$. Totik \cite{Totik2009,Totik2011} proved:

\begin{thm}Let $\mu $ have compact support $J$ and be regular. Let $I$ be a subinterval of $J$ satisfying the local Szeg\H{o} condition
\begin{gather*}
\int_{I}\big\vert \log \mu ^{\prime }(t) \big\vert dt<\infty.
\end{gather*}
 Then for a.e.\ $x\in I$, and all real $u$, $v$,
\begin{gather*}
\lim_{n\rightarrow \infty }\frac{\tilde{K}_{n}\left( \mu,x+\frac{u}{n\omega(x) },x+\frac{v}{n\omega (x) }\right) }{\tilde{K}_{n}(\mu,x,x) }=\mathbb{S}(u-v).
\end{gather*}
Here as above, $\omega $ is the equilibrium density of $J$.
\end{thm}

Totik actually showed that the asymptotic holds at any given $x$ which is a~Lebesgue point of both the measure $\mu $, and its local Szeg\H{o} function. This remains the most general result on explicit criteria for pointwise universality for compactly supported measures.

One question is whether a global condition such as regularity is needed, even if it is a weak one. Moreover, is there another way to treat general supports without using the polynomial pullback method of Totik? In \cite{Lubinsky2008}, a~method was introduced that avoids this. It uses basic tools of complex analysis and complex approximation, such as normal families, together with some of the theory of entire functions, and reproducing kernels.

Perhaps the most fundamental idea in this approach is the notion that since $K_{n}$ is a reproducing kernel for polynomials of degree $\leq n-1$, any scaled asymptotic limit of it must also be a reproducing kernel for a~suitable space. It turns out that when doing scaling in the bulk, the correct limit setting is Paley--Wiener space. This is the Hilbert space of entire functions~$g$ of exponential type at most $\pi $ whose restriction to the real line is in $L_{2}(\mathbb{R}) $. Here the sinc kernel is the reproducing kernel \cite[p.~95]{Stenger1993}:
\begin{gather*}
g(x) =\int_{-\infty }^{\infty }g(t) \mathbb{S}(x-t) dt,\qquad x\in \mathbb{R}.
\end{gather*}
This is the deeper reason for the appearance of the sinc kernel above. By using this idea, and complex analytic techniques, the author proved \cite{Lubinsky2008}:

\begin{thm}\label{thm2.5} Let $\mu $ have compact support $J$. Let $I$ be compact, and $\mu $ be absolutely continuous in an open set containing $I$. Assume that $\mu ^{\prime }$ is positive and continuous at each point of~$I$. The following are equivalent:
\begin{itemize}\itemsep=0pt
\item[$(I)$] Uniformly for $x\in I$ and $u$ in compact subsets of the real line,
\begin{gather}\label{eq2.3}
\lim_{n\rightarrow \infty }\frac{K_{n}\left( \mu,x+\frac{u}{n},x+\frac{u}{n}\right) }{K_{n}(\mu,x,x) }=1.
\end{gather}
\item[$(II)$] Uniformly for $x\in I$ and $u$, $v$ in compact subsets of the complex plane, we have
\begin{gather*}
\lim_{n\rightarrow \infty }\frac{K_{n}\left( \mu,x+\frac{u}{\widetilde{K}_{n}(\mu,x,x) },x+\frac{v}{\widetilde{K}_{n}(\mu,x,x) }\right) }{K_{n}(\mu,x,x) }=\mathbb{S}(u-v).
\end{gather*}
\end{itemize}
\end{thm}

One can weaken the condition of continuity of $\mu ^{\prime }$ to upper and lower bounds and then require~$x$ to be a Lebesgue point of~$\mu $, that is, we assume only
\begin{gather*}
\lim_{h,k\rightarrow 0+}\frac{\mu ([ x-h,x+k]) }{k+h}=\mu ^{\prime }(x).
\end{gather*}
The clear advantage of the theorem is that there is no global restriction on~$\mu$. The downside is that we still have to establish the ratio asymptotic~\eqref{eq2.3} for the Christof\/fel functions/ reproducing kernels, and to date, these have only been established in the stronger form of asymptotics for $K_{n}(\mu,x,x) $ itself. In the course of other work, Avila, Last and Simon \cite{Avilaetal2010,Simon2011} showed how to weaken the hypotheses on bounds of entire functions in~\cite{Lubinsky2008}.

With much more ef\/fort, and in particular a new uniqueness theorem for the sinc kernel, this set of methods also yields \cite{Lubinsky2012}: this is the only result that handles arbitrary measures with compact support.

\begin{thm}\label{thm2.6} Let $\mu $ have compact support, and let $J= \{ x\colon \mu ^{\prime }(x) >0 \} $. Let $\varepsilon >0$ and $r>0$. Then as $n\rightarrow \infty $,
\begin{gather*}
\operatorname{meas}\left\{ x\in J\colon \sup_{\vert u\vert,\vert v\vert
\leq r}\left\vert \frac{K_{n}\left( \mu,x+\frac{u}{\tilde{K}_{n} ( \mu,x,x) },x+\frac{v}{\tilde{K}_{n}(\mu,x,x) }\right) }{K_{n}(\mu,x,x) }-\mathbb{S}(u-v) \right\vert \geq \varepsilon \right\} \rightarrow 0.
\end{gather*}
\end{thm}

Here meas denotes linear Lebesgue measure, and in the supremum, $u$, $v$ are complex variables. Because convergence in measure implies convergence a.e.\ of subsequences, one obtains pointwise a.e.\ universality for subsequences, without any local or global assumptions on~$\mu $.

Another development involves pointwise universality in the mean \cite{Lubinsky2012B}, under some local conditions. Like Theorem~\ref{thm2.6}, the essential feature is the lack of global regularity assumptions:

\begin{thm}\label{thm2.7} Let $\mu $ have compact support. Assume that $I$ is an open interval in which for some $C>0$, $\mu ^{\prime }\geq C$ a.e.\ in $I$. Let $x\in I$ be a Lebesgue point of~$\mu $. Then for each $r>0$,
\begin{gather*}
\lim_{m\rightarrow \infty }\frac{1}{m}\sum_{n=1}^{m}\sup_{\vert u \vert,\vert v\vert \leq r}\left\vert \frac{K_{n}\left(
\mu,x+\frac{u}{\tilde{K}_{n}(\mu,x,x) },x+\frac{v}{\tilde{K}_{n}(\mu,x,x) }\right) }{K_{n}(\mu,x,x) }-\mathbb{S}(u-v) \right\vert =0.
\end{gather*}
In particular, this holds for a.e.\ $x\in I$.
\end{thm}

An obvious and important question is whether one can replace the convergence in measure in Theorem~\ref{thm2.6} with convergence a.e.\ or equivalently if one
really needs to take means as in Theorem~\ref{thm2.7}. Accordingly, we pose:

\begin{problem} Let $\mu $ have compact support. Assume that the support contains a non-empty interval $I$, in which $\mu ^{\prime }>0$ a.e. Is it true that for a.e.\ $x\in I$, and $u,v\in \mathbb{R}$, we have
\begin{gather*}
\lim_{n\rightarrow \infty }\frac{K_{n}\left( \mu,x+\frac{u}{\widetilde{K}_{n}(\mu,x,x) },x+\frac{v}{\widetilde{K}_{n} ( \mu,x,x) }\right) }{K_{n}(\mu,x,x) }=\mathbb{S}(u-v) ?
\end{gather*}
\end{problem}

My guess is that the answer is no. Vili Totik noted that the problem is open even for two basic situations:
\begin{itemize}\itemsep=0pt
\item[(i)] The restriction of $\mu $ to $I$ is Lebesgue measure, but there is no global assumption on $\mu $.
\item[(ii)] The support of $\mu $ is $[-1,1] $, and for example, $\mu $ is regular, but we do not assume a local Szeg\H{o} condition in~$I$.
\end{itemize}

Pointwise universality at a given point $x$ seems to usually require at least something like~$\mu ^{\prime }$ being continuous at $x$, or $x$ being a~Lebesgue point of~$\mu $. So it is quite a surprise that a~purely singularly continuous measure can exhibit this type of universality. This was shown by Jonathan Breuer~\cite{Breuer2011}, by starting with the Chebyshev weight of the second kind, $w(x) =\sqrt{1-x^{2}}$ on $[-1,1] $, and creating a new measure by sparsely perturbing the three term recurrence relation~\eqref{eq1.6} satisf\/ied by $\{ p_{n}(w,x) \} $. This special even weight has the recurrence relation
\begin{gather*}
xp_{n}(w,x) =\frac{1}{2}p_{n+1}(w,x)+\frac{1}{2}p_{n-1}(w,x).
\end{gather*}
For the perturbed measure $\mu $, we keep all $a_{n}=\frac{1}{2}$, and most $b_{n}=0$, but for a very sparse set of integers $\{ N_{j}\} $, set $b_{N_{j}}=v_{j}$, where $\{ v_{j}\} $ has limit~$0$:

\begin{thm} There exists a measure $\mu $ that is purely singularly continuous in $[-1,1] $, and has only mass points outside $[-1,1] $, such that for every $x\in (-1,1) $, and $u,v\in \mathbb{R}$,
\begin{gather*}
\lim_{n\rightarrow \infty }\frac{\pi \sqrt{1-x^{2}}}{n}K_{n}\left( \mu,x+\frac{u\pi \sqrt{1-x^{2}}}{n},x+\frac{v\pi \sqrt{1-x^{2}}}{n}\right) =\mathbb{S}(u-v).
\end{gather*}
\end{thm}

Another (and older) surprise is that universality can hold for measures supported on a Cantor set. This is a consequence of a result of Avila, Last and Simon on ergodic Jacobi matrices \cite{Avilaetal2010}. Let $\Omega $ be a compact metric space, $d\eta $ be a probability measure on $\Omega $, and $S\colon \Omega \rightarrow \Omega $ be an ergodic invertible map of $\Omega $ to itself. Let $A$, $B$ be continuous real valued functions on $\Omega $ with $\inf_{\Omega }A>0$. For each $\tau \in \Omega $, def\/ine a~Jacobi matrix
\begin{gather*}
J_{\tau }=\left[
\begin{matrix}
b_{1}(\tau) & a_{1}(\tau) & 0 & \cdots \\
a_{1}(\tau) & b_{2}(\tau) & a_{2} ( \tau) & \cdots \\
0 & a_{2}(\tau) & b_{3}(\tau) & \cdots \\
\vdots & \vdots & \vdots & \ddots
\end{matrix}
\right]
\end{gather*}
by
\begin{gather*}
a_{n}(\tau) =A\big( S^{n-1}\tau \big),\qquad b_{n}( \tau) =B\big( S^{n-1}\tau\big),\qquad n\geq 1.
\end{gather*}
Let $\mu _{\tau }$ be the spectral measure of $J_{\tau }$.

\begin{thm} Let $\{ J_{\tau }\} _{\tau \in S}$ be an ergodic Jacobi family as described above. Let $\Sigma _{ac}$ denote the common essential support of the a.c.\ spectrum of $J_{\tau }$, of positive Lebesgue measure. Then for a.e.\ $\tau \in S$ and for a.e.\ $x_{0}\in \sum_{ac}$, there is the universality limit
\begin{gather*}
\lim_{n\rightarrow \infty }\frac{K_{n}\left( \mu _{\tau },x_{0}+\frac{u}{n}
,x_{0}+\frac{v}{n}\right) }{n}=\frac{\rho (x_{0}) }{w_{\tau}(x_{0}) }\mathbb{S}( \rho (x_{0}) (u-v)),
\end{gather*}
with appropriate $\rho (x_{0}) $ and $w_{\tau}(x_{0}) $.
\end{thm}

One example is the almost Matthieu equation: let $\alpha $ be a f\/ixed irrational number, and
\begin{gather*}
a_{n}=1, \qquad b_{n}=2\lambda \cos ( \pi \alpha n+\theta),
\end{gather*}
$\lambda \in (-1,1) \backslash \{0\} $, $\Omega $ be the unit circle, $\tau =e^{i\theta }$, $S(\tau) =S(e^{i\theta }) =e^{i\theta }e^{i\pi \alpha }$ and $d\eta (\theta) =\frac{d\theta }{2\pi }$. The spectrum is purely absolutely continuous and a Cantor set. As far as the author is aware, this is the only known example where universality has been established for measures whose absolutely continuous spectrum is a Cantor set.

Another interesting phenomenon occurs at a jump discontinuity of a weight, where the universality is quite dif\/ferent from the sinc kernel. The precise form was obtained by Foulqui\'{e} Moreno, Mart\'{\i}nez-Finkelshtein, and Sousa~\cite{FoulquieMorenoetal2011}, using Riemann--Hilbert techniques. Let
\begin{gather*}
w_{c}(x) =h(x) (1-x) ^{\alpha }(1+x) ^{\beta }\Xi _{c}(x),\qquad x\in (-1,1),
\end{gather*}
where $h$ is positive on $[-1,1] $ and the restriction to $[-1,1] $ of a function analytic in a neighborhood of $[-1,1]$, while
\begin{gather*}
\Xi _{c}(x) =
\begin{cases}
1, & x\in \lbrack -1,0), \\
c^{2}, & x\in [0,1].
\end{cases}
\end{gather*}
Def\/ine
\begin{gather*}
G( a;z) =e^{-z/2}\sum_{k=0}^{\infty }\frac{(a) _{k}}{( k!) ^{2}}z^{k},
\end{gather*}
where $(a) _{k}=a( a+1) \cdots ( a+k-1) $ is the usual Pochhammer symbol. Def\/ine $\lambda =\frac{i}{\pi }\log c$ and if $x\neq y$, the kernel
\begin{gather*}
\mathbb{K}(x,y) =\frac{1}{i\pi h(0) }\frac{\log c}{c^{2}-1}\frac{G( 1+\lambda ;2\pi ix) G( \lambda ;2\pi iy) -G( \lambda ;2\pi ix) G( 1+\lambda ;2\pi iy)}{x-y},
\end{gather*}
while
\begin{gather*}
\mathbb{K}( x,x) =\frac{2}{h(0) }\frac{\log c}{c^{2}-1}\big[ G^{\prime }( 1+\lambda ;2\pi ix) G( \lambda ;2\pi ix) -G^{\prime }( \lambda ;2\pi ix) G(1+\lambda ;2\pi ix) \big].
\end{gather*}
Foulqui\'{e} Moreno, Mart\'{\i}nez-Finkelshtein, and Sousa established detailed asymptotics for ortho\-go\-nal polynomials for the weight $w_{c}$, and hence deduced:

\begin{thm} Let $c>0$, $c\neq 1$, and $\delta \in (0,1) $. Then uniformly for $u,v\in ( -\delta,\delta) $,
\begin{gather*}
\lim_{n\rightarrow \infty }\frac{\pi }{n}K_{n}\left( w_{c},\frac{\pi u}{n},\frac{\pi v}{n}\right) =\mathbb{K}(u,v).
\end{gather*}
\end{thm}

Note that when we let $c\rightarrow 1$, we recover the sine kernel. In \cite{Lubinsky2009B}, the author attempted to classify possible limiting kernels arising from universality limits for a broad class of compactly supported weights. If $\mu $ has compact support, and~$O$ is some open set inside the support in which $\mu $ is absolutely continuous, while $\mu ^{\prime }$ is bounded above and below by positive constants there, then I proved that every function $F(u,v) $ that is a scaled limit of some subsequence of $\left\{ \frac{K_{n}\left( \mu,x+\frac{u}{\widetilde{K}_{n}(\mu,x,x)},x+\frac{v}{\widetilde{K}_{n}(\mu,x,x) }\right) }{K_{n}(\mu,x,x) }\right\} $ is the reproducing kernel of some de Branges space that is equivalent to a~Paley--Wiener space. Yes this is technical, see~\cite{Lubinsky2009B} for details. There is numerical evidence that the new kernel obtained by Foulqui\'{e} Moreno, Mart\'{\i}nez-Finkelshtein, and Sousa f\/its this mold. There is a also a natural connection between their results and those on Toeplitz determinants for weights with more general Fisher--Hartwig type singularities.

The result of Breuer on singularly continuous measures suggests that universality limits can be preserved when we perturb measures. In \cite{Breueretal2014}, Breuer, Last, and Simon showed that when we start with a~base measure $\mu $ for which we have bulk universality, and moderately perturb the recurrence coef\/f\/icients of~$\mu$ to obtain a new measure $\tilde{\mu}$, then universality persists, and moreover, for random perturbations, universality persists under weaker conditions. Here is one of their deterministic results:

\begin{thm}Let $\mu $ be a compactly supported measure with recurrence coefficients $\{a_{n}\}$, $\{b_{n}\}$. Let $w=\mu ^{\prime }$ denote its Radon--Nikodym derivative. Assume that $A\subset \mathbb{R}$ and that for a.e.\ $x_{0}\in A$, there exists a number $\rho (x_{0})$ such that
\begin{gather}\label{eq2.4}
\lim_{n\rightarrow \infty }\frac{K_{n}\left( \mu,x_{0}+\frac{u}{n},x_{0}+\frac{v}{n}\right) }{n}=\frac{\rho (x_{0}) }{w(x_{0}) }\mathbb{S}( \rho (x_{0}) (u-v)),
\end{gather}
uniformly for $u$, $v$ in compact subsets of $\mathbb{C}$. Let $\{ \beta _{k}\} \subset \mathbb{R}$ satisfy
\begin{gather*}
\sum_{k=1}^{\infty }\vert \beta _{k}\vert <\infty.
\end{gather*}
Let $\tilde{\mu}$ be the measure with recurrence coefficients $\{a_{n}\}$, $\{ b_{n}+\beta _{n}\} $. Then \eqref{eq2.4} is also true for the measure $\tilde{\mu}$, for a.e.\ $x_{0}\in A$, and $u$, $v$ in compact subsets of~$\mathbb{C}$, but with $\rho (x_{0}) $ and $w(x_{0}) $ replaced by appropriately modified $\tilde{\rho}(x_{0}) $ and $\tilde{w}(x_{0}) $.
\end{thm}

In many of the proofs, asymptotics of Christof\/fel functions $1/K_{n}(\mu,x,x) $ play a key role. Those asymptotics are often established with the aid of the variational principle
\begin{gather*}
K_{n}(\mu,x,x) =\sup_{\deg (P) \leq n-1}\frac{P(x) ^{2}}{\int P(t) ^{2}d\mu (t) },
\end{gather*}
which immediately implies that $K_{n}(\mu,x,x) $ is monotone decreasing in $\mu $. In \cite[Theorem~1.1, p.~111]{Lubinsky2013}, the author established a~similar variational principle for the general $m$-point correlation function $R_{m,n}$. Its formulation involves $\mathcal{AL}_{n}^{m}$, the alternating polynomials of degree at most~$n$ in~$m$ variables. We say that $P\in \mathcal{AL}_{n}^{m}$ if
\begin{gather*}
P\left( x_{1},x_{2},\dots,x_{m}\right) =\sum_{0\leq j_{1},j_{2},\dots,j_{m}\leq n}c_{j_{1}j_{2}\cdots j_{m}}x_{1}^{j_{1}}x_{2}^{j_{2}}\cdots x_{m}^{j_{m}},
\end{gather*}
so that $P$ is a polynomial of degree $\leq n$ in each of its $m$ variables, and in addition is \textit{alternating}, so that for every pair $(i,j)$ with $1\leq i<j\leq m$,
\begin{gather*}
P( x_{1},\dots,x_{i},\dots,x_{j},\dots,x_{m}) =-P(x_{1},\dots,x_{j},\dots,x_{i},\dots,x_{m}).
\end{gather*}
Thus swapping variables changes the sign.

\begin{thm}\label{thm2.13}
\begin{gather*}
\det [ K_{n}( \mu,x_{i},x_{j})] _{1\leq i,j\leq m}=m!\sup_{P\in \mathcal{AL}_{n-1}^{m}}\frac{( P(x_{1},x_{2},\dots,x_{m})) ^{2}}{\int ( P(t_{1},t_{2},\dots,t_{m})) ^{2}d\mu ( t_{1}) d\mu ( t_{2}) \cdots d\mu ( t_{m}) }.
\end{gather*}
\end{thm}

An immediate consequence is monotonicity of the \textit{unweighted} $m$-point correlation function for measures:

\begin{coro}
\begin{gather*}
\hat{R}_{m,n} ( \mu ;x_{1},x_{2},\dots,x_{m})\\
\qquad{} :=\frac{n!}{\left( n-m\right) !}\frac{1}{Z_{n}}\int \cdots \int \left(\prod\limits_{1\leq i<j<n}( x_{i}-x_{j}) ^{2}\right) d\mu (x_{m+1}) d\mu ( x_{m+2}) \cdots d\mu ( x_{n}) \\
\qquad{} =\det [ K_{n}( \mu,x_{i},x_{j}) ] _{1\leq i,j\leq m}
\end{gather*}
is a monotone decreasing function of $\mu $.
\end{coro}

Note that when $\mu $ is absolutely continuous and $d\mu (x) =w(x) dx$, then
\begin{gather*}
R_{m,n}\left( w;x_{1},x_{2},\dots,x_{m}\right) =\hat{R}_{m,n}\left( \mu
;x_{1},x_{2},\dots,x_{m}\right) w\left( x_{1}\right) w\left( x_{2}\right)
\cdots w\left( x_{m}\right).
\end{gather*}

The proof of Theorem~\ref{thm2.13} is based on multivariate alternating orthogonal polynomials built from~$\mu $. As the author found out after writing \cite{Lubinsky2013}, these polynomials were known before \cite{KarlinMcGregor1962}, \cite[p.~182, Theorem~3.8.6]{Simon2015}, though it seems that the monotonicity property is new. One consequence is one-sided universality without any restrictions on the measure \cite[Theorem~2.2, p.~116]{Lubinsky2013}:

\begin{thm}Let $\mu $ have compact support $\mathcal{K}$, and let $\omega $ denote the density of the equilibrium measure for $\mathcal{K}$ in the interior $\mathcal{K}^{o}$ of $\mathcal{K}$. Let $J=\{ x\colon \mu ^{\prime }(x)>0\}$. Let $m\geq 1$. For a.e.\ $x\in J\cap \mathcal{K}^{o}$, and for all real $u_{1},u_{2},\dots,u_{m}$,
\begin{gather*}
\liminf_{n\rightarrow \infty }\left( \frac{\mu ^{\prime }(x) }{n\omega (x) }\right) ^{m}R_{m,n}\left( \mu,x+\frac{u_{1}}{n\omega (x) },\dots,x+\frac{u_{m}}{n\omega (x) }\right) \geq \det ( \mathbb{S}( u_{i}-u_{j}) ) _{1\leq i,j\leq m}.
\end{gather*}
\end{thm}

We note that an analogous upper bound for limsup's is also proved in \cite{Lubinsky2013}, with $\omega $ replaced by a density formed from taking inf's of equilibrium densities of compact subsets $L\subset J$ such that $\mu _{|L}$ is regular.

A very recent and exciting development is the treatment of endpoint and interior power singularities for quite general measures (cf.\ \cite{Danka2016,Danka2016B,DankaTotik2016}). Tivadar Danka, a~student of Vili Totik, has developed this theory using a mix of techniques, including Riemann--Hilbert, and reproducing kernel ideas. The kernels are modif\/ied Bessel kernels. (In the case of general varying weights that include a factor of $\vert x\vert ^{\alpha }$, Kuijlaars and Vanlessen~\cite{KuijlaarsVanlessen2003} earlier obtained universality results, extending several earlier works.) Let
\begin{gather*}
\mathbb{L}_{\alpha }(u,v) =\frac{\sqrt{uv}}{2}\frac{J_{\frac{\alpha +1}{2}}(u) J_{\frac{\alpha -1}{2}}(v) -J_{\frac{\alpha +1}{2}}(v) J_{\frac{\alpha -1}{2}}(u) }{u-v}\qquad \text{if} \quad u,v\geq 0, \quad u\neq v
\end{gather*}
and
\begin{gather*}
\mathbb{L}_{\alpha }(u,v) =\mathbb{L}_{\alpha }(\vert u\vert,\vert v\vert) \qquad \text{if} \quad \min \{u,v\} <0,\quad u\neq v.
\end{gather*}
Also, along the diagonal, let
\begin{gather*}
\mathbb{L}_{\alpha }( u,u) =\frac{\vert u\vert }{2}\left( J_{\frac{\alpha +1}{2}}^{\prime }(u) J_{\frac{\alpha -1}{2}}(u) -J_{\frac{\alpha +1}{2}}(u) J_{\frac{\alpha -1}{2}}^{\prime }(u) \right).
\end{gather*}
Since these kernels can include non-integer powers of $z$, def\/ine also their entire cousins
\begin{gather*}
\mathbb{L}_{\alpha }^{\ast }(u,v) =\frac{\mathbb{L}_{\alpha}(u,v) }{u^{\alpha /2}v^{\alpha /2}}.
\end{gather*}
Following is Danka's result for interior power singularities \cite[Theorem~1.3]{Danka2016}:

\begin{thm}Let $\mu $ have compact support $J$ and be regular. Let $x_{0}$ be an interior point of~$J$, $\delta>0$, $\alpha >-1$, and suppose that in $(x_{0}-\delta,x_{0}+\delta)$, $\mu $ is absolutely continuous, with
\begin{gather*}
\mu ^{\prime }(x) =w(x) \vert x-x_{0}\vert ^{\alpha }, \qquad x\in ( x_{0}-\delta,x_{0}+\delta),
\end{gather*}
where $w$ is strictly positive and continuous at $x_{0}$. Let $\omega $ denote the equilibrium density of~$J$. Then
\begin{gather*}
\lim_{n\rightarrow \infty }\frac{K_{n}\left( \mu,x_{0}+\frac{u}{n},x_{0}+\frac{v}{n}\right) }{K_{n}( x_{0},x_{0}) }=\frac{\mathbb{L}_{\alpha }^{\ast }( \pi \omega (x_{0}) u,\pi \omega(x_{0}) v) }{\mathbb{L}_{\alpha }^{\ast }(0,0) }.
\end{gather*}
\end{thm}

For endpoint singularities, Danka proved \cite[Theorem~1.4]{Danka2016}

\begin{thm}Let $\mu $ have compact support $J$ and be regular. Let $x_{0}$ be a right endpoint of $J$, in the sense that $J\cap ( x_{0},x_{0}+\varepsilon) =\varnothing $ for some $\varepsilon >0$. Assume that $\delta >0$, $\alpha >-1$, and suppose that in $(x_{0}-\delta,x_{0}]$, $\mu $ is absolutely continuous, with
\begin{gather*}
\mu ^{\prime }(x) =w(x) \vert x-x_{0}\vert ^{\alpha }, \qquad x\in (x_{0}-\delta,x_{0}],
\end{gather*}
where $w$ is strictly positive and left-continuous at $x_{0}$. Let $\omega $ denote the equilibrium density of~$J$. Then
\begin{gather*}
\lim_{n\rightarrow \infty }\frac{K_{n}\left( \mu,x_{0}-\frac{u}{2n^{2}},x_{0}-\frac{v}{2n^{2}}\right) }{K_{n}( x_{0},x_{0}) }=\frac{\mathbb{J}_{\alpha }^{\ast }\left( L^{2}u,L^{2}v\right) }{\mathbb{J}_{\alpha}^{\ast }(0,0) },
\end{gather*}
where
\begin{gather*}
L=\lim_{x\rightarrow x_{0}-}\sqrt{2}\pi \vert x-x_{0}\vert^{1/2}\omega (x).
\end{gather*}
\end{thm}

We began this section with one powerful illustration of the Deift--Zhou Riemann--Hilbert method, for Jacobi weights. We close this section with another, this time by Shuai-Xia Xu, Yu-Qiu Zhao, and Jian-Rong Zhou~\cite{Xuetal2011}, for the weights
\begin{gather}\label{eq2.5}
w(x) =\exp \big( {-}\big( 1-x^{2}\big) ^{-\Delta }\big), \qquad x\in (-1,1).
\end{gather}
For these weights, universality in the bulk follows from, for example, Theorem~\ref{thm2.2}. However, the authors obtained rates, but the really interesting features occur near $\pm 1$. Because $w$ decays rapidly near~$\pm 1$, we need to use the Mhaskar--Rakhmanov--Saf\/f interval $[-\beta _{n},\beta_{n}]$, where $\beta _{n}$ is the root of
\begin{gather*}
\int_{0}^{\beta _{n}}\frac{x^{2}}{\left( 1-x^{2}\right) ^{\Delta +1}\sqrt{\beta _{n}^{2}-x^{2}}}dx=\frac{n\pi }{2\Delta }.
\end{gather*}
(See Section~\ref{section6} for more on Mhaskar--Rakhmanov--Saf\/f intervals.) There is a~complete asymptotic expansion for $\beta _{n}$, the f\/irst terms of which are
\begin{gather*}
\beta _{n}=1-\left[ \frac{1}{2}\left( \frac{\Gamma \left( \Delta +\frac{1}{2}\right) }{\sqrt{\pi }\Gamma ( \Delta) }\right) ^{\frac{1}{\Delta
+\frac{1}{2}}}\right] n^{-\frac{1}{\Delta +\frac{1}{2}}}( 1+o(1)).
\end{gather*}

\begin{thm}Let $\Delta >0$ and $w$ be given by \eqref{eq2.5}. If $\Delta \neq \frac{1}{2}$, let
\begin{gather*}
\varepsilon _{n}=n^{-\min \left\{ 1,\frac{1}{\Delta +\frac{1}{2}}\right\} }
\end{gather*}
while if $\Delta =\frac{1}{2}$, let $\varepsilon_{n}=n^{-1}\log n$.
\begin{itemize}\itemsep=0pt
\item[$(a)$] For $x\in (-1,1) $, and uniformly for $u$, $v$ in compact subsets of $\mathbb{R}$,
\begin{gather*}
\frac{\pi \sqrt{1-x^{2}}}{n}K_{n}\left( w,x+\frac{u\pi \sqrt{1-x^{2}}}{n},x+\frac{v\pi \sqrt{1-x^{2}}}{n}\right) =\mathbb{S}(u-v) +O(\varepsilon _{n}).
\end{gather*}
\item[$(b)$] There is a number $B_{0}$ such that uniformly for $u$, $v$ in compact subsets of $\mathbb{R}$,
\begin{gather*}
B_{0}n^{-\frac{4\Delta +6}{6\Delta +3}}K_{n}\left( w,\beta _{n}+uB_{0}n^{-\frac{4\Delta +6}{6\Delta +3}},\beta _{n}+vB_{0}n^{-\frac{4\Delta +6}{6\Delta +3}}\right) =\mathbb{A}(u,v) +O\left( n^{-\frac{4\Delta}{6\Delta +3}}\right).
\end{gather*}
\end{itemize}
\end{thm}

We note that $B_{0}$ is explicitly given. The interesting feature is that we obtain the Airy, rather than Bessel, kernel, but as the edge $\beta _{n}$ is
``soft'', this should not be surprising. We close this section with some open problems:

\begin{problem}Let $\mu $ be a measure with compact support.
\begin{itemize}\itemsep=0pt
\item[$(a)$] Assume that $I$ is a subinterval of the support and that $\mu ^{\prime }$ is bounded above and below by positive constants in~$I$. Describe the set of limiting kernels arising from scaling limits around points in~$I$.
\item[$(b)$] More generally, without any restrictions on~$\mu $, describe the set of limiting kernels arising from scaling limits around points in $\operatorname{supp}[\mu]$.
\end{itemize}
\end{problem}

As we noted above, something like de Branges spaces might play a~role~\cite{Lubinsky2009B}. Following is a~more specif\/ic problem:

\begin{problem}
Investigate the universality limits for a compactly supported measure $\mu $ around a ``strong'' interior zero, such as
\begin{gather*}
\mu ^{\prime }(x) =\exp \left( -\exp \left( \cdots \exp \left(\vert x\vert ^{-\alpha }\right) \right) \right), \qquad x\in (-1,1), \quad \alpha >0.
\end{gather*}
\end{problem}

\section{Applications in orthogonal polynomials}\label{section3}

Since $\mathbb{S}(z) =\frac{\sin \pi z}{\pi z}$ has zeros at all the non-zero integers, the universality limit~\eqref{eq1.17} suggests that these should attract zeros of the scaled reproducing kernels, which include zeros of orthogonal polynomials. One can make this rigorous using Hurwitz' theorem if the universality limit holds uniformly for $u$, $v$ in compact subsets of the plane. Alternatively, if one only knows this only for real~$u$,~$v$, we can use the intermediate value theorem and more elementary techniques. This connection has been explored by several authors. In some sense it goes back to Freud~\cite{Freud1971}, though it was Eli Levin who discovered and formalized the idea in~\cite{LevinLubinsky2008B}. We present only the most general known result, due to V.~Totik \cite[Theorem~2.1]{Totik2009}. Other aspects of ``clock spacing'' of zeros, relating to universality, have been explored by Barry Simon and his collaborators in \cite{Avilaetal2010,Simon2005,Simon2008,Simon2011}. The zeros of $p_{n}(\mu,x)$ are denoted by
\begin{gather*}
-\infty <x_{nn}<x_{n,n-1}<\cdots <x_{1n}<\infty.
\end{gather*}

\begin{thm}\label{thm3.1}Let $\mu $ be a regular measure with compact support $J\subset \mathbb{R}$. Let $\mathcal{K}$ be a compact subset of the interior of $J$. Assume that $\mu $ is absolutely continuous in an open set containing $\mathcal{K}$, and that $\mu ^{\prime }$ is positive and continuous at each point of $\mathcal{K}$. Let $\omega $ denote the equilibrium density of~$J$. Let $L\geq 1$. Then
\begin{gather}\label{eq3.1}
\lim_{n\rightarrow \infty }n( x_{kn}-x_{k+1,n}) \omega (x) =1
\end{gather}
uniformly in $x\in \mathcal{K}$ and $\vert x_{kn}-x\vert \leq L/n$.
\end{thm}

Totik~\cite[Theorem~2.3]{Totik2009} also proved a version of this result where continuity of $\mu ^{\prime }$ is replaced by integrability of $\log \mu ^{\prime }$ in a neighborhood of a point $x_{0}$, which is assumed to be a Lebesgue point of appropriate functions. For weights with jump discontinuities, spacing of zeros was explored in~\cite{FoulquieMorenoetal2011}, where it is shown that there is not ``clock behavior'' around the discontinuity, so that~\eqref{eq3.1} fails. For weights with interior power singularities, Danka has completed the task~\cite{Danka2016B}.

Our next result \cite[Theorem~1.2, p.~72]{LevinLubinsky2008B} concerns asymptotics of the zeros close to $1$, which follows from Theorem~\ref{thm2.3}. We
denote the positive zeros of the Bessel function $J_{\alpha }$ by
\begin{gather*}
0<j_{\alpha,1}<j_{\alpha,2}<j_{\alpha,3}<\cdots.
\end{gather*}

\begin{thm}\label{thm3.2}
Let $\mu $ be a finite positive Borel measure on $(-1,1)$ that is regular. Assume that for some $\rho >0$, $\mu $ is absolutely continuous in $\mathcal{K}=[ 1-\rho,1] $, and in $\mathcal{K}$, its absolutely continuous component has the form $w=hw_{\alpha,\beta }$, where $\alpha,\beta >-1$. Assume that $h(1)>0$ and~$h$ is continuous at~$1$. Then for each fixed $k\geq 1$,
\begin{gather*}
\lim_{n\rightarrow \infty }n\sqrt{1-x_{kn}^{2}}=j_{\alpha,k}
\end{gather*}
and
\begin{gather*}
\lim_{n\rightarrow \infty }n^{2}( x_{kn}-x_{k+1,n}) =\frac{1}{2}\big( j_{\alpha,k+1}^{2}-j_{\alpha,k}^{2}\big).
\end{gather*}
\end{thm}

Analogues of Theorems~\ref{thm3.1} and \ref{thm3.2} for exponential and varying weights, based on universality limits, have been explored in~\cite{Deiftetal1999B,LevinLubinsky2008B}. In particular, Deift et al.\ give remarkably precise results for the largest and smallest zeros and for spacing of interior zeros.

It is noteworthy that asymptotics of zeros are in some sense equivalent to universality. Another equivalence condition involves the dif\/ferentiated kernel
\begin{gather*}
K_{n}^{( r,s) }(\mu,x,y) =\sum_{k=0}^{n-1}p_{k}^{(r)}(x) p_{k}^{(s) }(y),
\end{gather*}
and its normalized cousin
\begin{gather*}
\tilde{K}_{n}^{(r,s) }(\mu,x,y) =\mu ^{\prime
}(x) ^{1/2}\mu ^{\prime }(y) ^{1/2}K_{n}^{(r,s) }(x,y)
\end{gather*}
for non-negative integers $r$, $s$. We also def\/ine
\begin{gather*}
\tau _{r,s}= \begin{cases}
0, & \text{$r+s$ odd}, \\
\frac{( -1) ^{( r-s) /2}}{r+s+1}, &
\text{$r+s$ even}.
\end{cases}
\end{gather*}

Given a f\/ixed real number $\xi $, we let $\{ t_{j,n}\} _{j}=\{ t_{j,n}( \xi)\} _{j}$ denote the $n-1$ or $n$ zeros of the ``quasi-orthogonal'' polynomial
\begin{gather*}
p_{n}(\mu,t) p_{n-1}( \mu,\xi) -p_{n-1}( \mu,t) p_{n}( \mu,\xi ).
\end{gather*}
They are real, simple, and interlace the zeros of $p_{n}(\mu,x)$ \cite[p.~19]{Freud1971}. We assume they are ordered so that
\begin{gather*}
\cdots <t_{-2,n}(\xi) <t_{-1,n}(\xi) <t_{0,n}(\xi) =\xi <t_{1,n}(\xi) <t_{2,n}(\xi) <\cdots.
\end{gather*}
Of course it is possible that all $t_{k,n}$, other than $\xi $, lie to the left or right of $\xi $. Below $\operatorname{dist}(x,J) $ denotes the distance from a point $x$ to a set $J$. The following equivalences appeared in \cite[Theorem~1.3, p.~181]{LevinLubinsky2010B}:

\begin{thm}Let $\mu $ be a finite positive Borel measure on the real line with compact support. Let $J\subset \operatorname{supp}[ \mu]$ be compact, and such that $\mu $ is absolutely continuous in an open set containing~$J$. Assume that $\mu ^{\prime }$ is positive and continuous at each point of $J$. The following are equivalent:
\begin{itemize}\itemsep=0pt
\item[$(I)$] Uniformly for $x\in J$ and $u$, $v$ in compact subsets of the complex plane, we have
\begin{gather*}
\lim_{n\rightarrow \infty }\frac{K_{n}\left( \mu,x+\frac{u}{\widetilde{K}_{n}(\mu,x,x) },x+\frac{v}{\widetilde{K}_{n}( \mu,x,x) }\right) }{K_{n}(\mu,x,x) }=\mathbb{S}(u-v).
\end{gather*}
\item[$(II)$] For each $r,s\geq 0$, and uniformly for $x\in J$,
\begin{gather*}
\lim_{n\rightarrow \infty }\frac{\tilde{K}_{n}^{(r,s) }(\mu,x,x) }{\widetilde{K}_{n}(\mu,x,x) ^{r+s+1}}=\pi^{r+s}\tau _{r,s}.
\end{gather*}
\item[$(III)$] There exists $L>0$, such that uniformly in $n$ and for $x $ with $\operatorname{dist}(x,J) \leq L/n$,
\begin{gather*}
K_{n}^{(1,0) }(\mu,x,x) =o\big( n^{2}\big).
\end{gather*}
\item[$(IV)$] For each fixed $j$, we have uniformly in $x\in J$,
\begin{gather*}
\lim_{n\rightarrow \infty }( t_{n,j+1}(x) -t_{n,j}(x)) \tilde{K}_{n}(\mu,x,x) =1.
\end{gather*}
\end{itemize}
\end{thm}

Remarkably, asymptotics for $K_{n}^{(r,s) }$ for $r,s=0,1$, are just what are needed for studying zeros of random polynomials, expressed as linear combinations of $\{ p_{j}(\mu,x)\} $. This was observed by Igor Pritsker, and used in \cite{Lubinskyetal2016,PritskerXie2016}.

\section{Bulk limits and sine kernels for other orthogonal systems}\label{section4}

The bulk limit, in the form \eqref{eq1.17}, has been extended to other systems of orthogonal ``polyno\-mials''. Of course in many of these cases, the connection to random matrices is lost. We brief\/ly cover these in this section:

\subsection{Unit circle}\label{section4.1}

Let $\mu $ be a f\/inite positive Borel measure on $[-\pi,\pi )$ with inf\/initely many points in its support. Then we may def\/ine orthonormal polynomials $\phi _{n}(z) =\kappa _{n}z^{n}+\cdots$, $\kappa_{n}>0$, $n=0,1,2,\dots $ satisfying the orthonormality conditions
\begin{gather}\label{eq4.1}
\frac{1}{2\pi }\int_{-\pi }^{\pi }\phi _{n}(z) \overline{\phi_{m}(z) }d\mu (\theta) =\delta _{mn},
\end{gather}
where $z=e^{i\theta }$. The normalization by $2\pi $ is part of the ``culture'' of this topic. Regularity in the sense of Stahl, Totik and Ullman takes the form
\begin{gather*}
\lim_{n\rightarrow \infty }\kappa _{n}^{1/n}=1.
\end{gather*}
The $n$th reproducing kernel for $\mu $ is
\begin{gather}\label{eq4.2}
K_{n}( \mu,z,u) =\sum_{j=0}^{n-1}\phi _{j}(z)\overline{\phi _{j}(u) },
\end{gather}
with normalized cousin
\begin{gather}\label{eq4.3}
\tilde{K}_{n}\big( \mu,e^{i\theta },e^{is}\big) =\mu ^{\prime }(\theta) ^{1/2}\mu ^{\prime }(s) ^{1/2}K_{n}\big( \mu,e^{i\theta },e^{is}\big).
\end{gather}
Sometimes, we identify $\mu ^{\prime }(\theta) $ with $\mu ^{\prime }( e^{i\theta }) $. Eli Levin and the author proved \cite{LevinLubinsky2007}:

\begin{thm}\label{thm4.1} Let $\mu $ be a finite positive Borel measure on $[-\pi,\pi )$ that is regular. Let $J\subset (-\pi,\pi )$ be compact, and such that $\mu $ is absolutely continuous in an open set containing $J$. Assume moreover, that $\mu ^{\prime }$ is positive and continuous at each point of~$J$. Then uniformly for $\theta \in J$ and $u$, $v$ in compact subsets of the plane, we have
\begin{gather*}
\lim_{n\rightarrow \infty }\frac{1}{n}\tilde{K}_{n}\left( \mu,e^{i\left(\theta +\frac{2\pi u}{n}\right) },e^{i\left( \theta +\frac{2\pi \bar{v}}{n}\right) }\right) =e^{i\pi (u-v) }\mathbb{S}(u-v).
\end{gather*}
Equivalently,
\begin{gather*}
\lim_{n\rightarrow \infty }\frac{K_{n}\left( \mu,z\left( 1+\frac{i2\pi u}{n}\right),z\left( 1+\frac{i2\pi \bar{v}}{n}\right) \right) }{K_{n}\left(\mu,z,z\right) }=e^{i\pi (u-v) }\mathbb{S}(u-v),
\end{gather*}
uniformly for $u$, $v$ in compact subsets of the complex plane and $z=e^{i\theta }$, $\theta \in J$.
\end{thm}

Of course, if $\mu ^{\prime }( -\pi ) =\mu ^{\prime }( \pi) $, and $\mu ^{\prime }$ is continuous at $\pm \pi $, then this result also holds at $z=e^{\pm i\pi }$.

\begin{coro} Let $r$, $s$ be non-negative integers and
\begin{gather*}
K_{n}^{(r,s) }(\mu,z,z) =\sum_{k=0}^{n-1}\phi _{k}^{(r) }(z) \overline{\phi _{k}^{(s) }(z) }.
\end{gather*}
Then uniformly for $\theta \in J$, $z=e^{i\theta }$,
\begin{gather*}
\lim_{n\rightarrow \infty }\frac{z^{r-s}}{n^{r+s}}\frac{K_{n}^{(r,s) }( \mu,z,z) }{K_{n}(\mu,z,z)}=\frac{1}{r+s+1}.
\end{gather*}
\end{coro}

\subsection{Arcs of the unit circle}\label{section4.2}

Let $\alpha \in \left( 0,\pi \right) $ and let our arc be $\Delta _{\alpha}=\{ e^{i\theta }\colon \theta \in [ \alpha,2\pi -\alpha]\} $. Let $\mu $ be a f\/inite positive Borel measure on $\Delta_{\alpha }$ (or equivalently on $[\alpha,2\pi -\alpha ] $) with inf\/initely many points in its support. Then we may def\/ine orthonormal polynomials $\phi _{n}(z) =\kappa _{n}z^{n}+\cdots$, $\kappa _{n}>0$, $n=0,1,2,\dots $ satisfying~\eqref{eq4.1} with $[ -\pi,\pi]$ replaced by $[ \alpha,2\pi -\alpha] $. The reproducing kernel and its normalized cousin are given by~\eqref{eq4.2},~\eqref{eq4.3}. For $\theta \in [\alpha,2\pi -\alpha]$, let
\begin{gather*}
T(\theta) =\frac{\sin \frac{\theta }{2}}{\sqrt{\cos ^{2}\frac{\alpha }{2}-\cos ^{2}\frac{\theta }{2}}}.
\end{gather*}
$T(\theta) /( 2\pi) $ is the density of the equilibrium measure for $\Delta _{\alpha }$. The author and Nguyen \cite{LubinskyNguyen2013} proved

\begin{thm} Let $\alpha \in ( 0,\pi) $, and let $\mu $ be a finite positive Borel measure on $[ \alpha,2\pi-\alpha ] $ that is regular. Let $J\subset (\alpha,2\pi -\alpha) $ be compact, and be such that $\mu $ is absolutely continuous in an open set containing~$J$. Assume moreover, that $\mu ^{\prime }$ is positive and continuous at each point of~$J$. Then uniformly for $\theta _{0}\in J$ and $u$, $v$ in compact subsets of the complex plane, we have
\begin{gather*}
\lim_{n\rightarrow \infty }\frac{K_{n}\left( \mu,e^{i\left( \theta _{0}+\frac{2\pi u}{n}\right) },e^{i\left( \theta _{0}+\frac{2\pi \bar{v}}{n}\right) }\right) }{K_{n}\left( \mu,e^{i\theta _{0}},e^{i\theta _{0}}\right)
}=e^{i\pi (u-v) }\mathbb{S}( (u-v) T(\theta _{0})).
\end{gather*}
\end{thm}

This reduces to Theorem~\ref{thm4.1} as $\alpha \rightarrow 0+$. Applications to asymptotics of zeros and dif\/fe\-ren\-tia\-ted kernels were also presented in~\cite{LubinskyNguyen2013}.

\subsection{Smooth closed contours}\label{section4.3}

Consider a smooth closed contour $\Gamma =\{ \gamma (s)\colon s\in [ 0,L]\} $, where $L>0$. $\Gamma $ is assumed to be ``smooth'' in the following sense: $\gamma^{\prime \prime }$ exists and is continuous on $[0,L] $, and satisf\/ies a Lipschitz condition of some positive order $\beta >0$. Thus, for some $C>0$,
\begin{gather*}
\big\vert \gamma ^{\prime \prime }(s) -\gamma ^{\prime \prime }(t)\big\vert \leq C\vert t-s\vert ^{\beta }, \qquad s,t\in [0,L].
\end{gather*}
In addition, we assume that $\gamma $ is periodic on $[0,L] $, so that $\gamma ^{(j) }(0) =\gamma ^{(j) }(L) $, $j=0,1,2$. In Suetin's 1966 terminology \cite{Suetin1966}, $\Gamma \in C( 2,\beta) $.

We denote the exterior of $\Gamma $ by $D$, and denote the conformal map of~$D$ onto the exterior of the unit ball by~$\Phi $, normalized by $\Phi(\infty) =\infty $, and $\Phi ^{\prime }(\infty) >0$. We assume that $\mu $ is a f\/inite positive Borel measure on $\Gamma $, and $\{ \phi _{n}\}$ are orthonormal polynomials for $\mu $, so that
\begin{gather*}
\frac{1}{2\pi }\int_{\Gamma }\phi _{n}(z) \overline{\phi_{m}(z) }d\mu (z) =\delta _{mn}.
\end{gather*}
As usual, $K_{n}$ is given by~\eqref{eq4.2}. The author and Levin proved \cite{LevinLubinsky2016}:

\begin{thm} Let $\Gamma $ be a simple closed curve in the complex plane, of class $C( 2,\beta) $, for some $\beta \in (0,1) $. Let $\mu $ be a finite positive Borel measure on $\Gamma $ that is regular. Let $\Gamma _{1}$ be a closed proper subarc of $\Gamma $, such that $\mu $ is absolutely continuous with respect to arclength, in an open arc containing $\Gamma _{1}$, and the Radon--Nikodym derivative $\mu^{\prime }$ $($with respect to arclength$)$ is positive and continuous in that open subarc. Then uniformly for $z_{0}\in \Gamma _{1}$ and $u$, $v$ in compact subsets of~$\mathbb{C}$,
\begin{gather*}
\lim_{n\rightarrow \infty }\frac{K_{n}\left( \mu,z_{0}+\frac{2\pi iu}{n}\frac{\Phi (z_{0}) }{\Phi ^{\prime }(z_{0}) },z_{0}+\frac{2\pi i\bar{v}}{n}\frac{\Phi (z_{0}) }{\Phi ^{\prime}(z_{0}) }\right) }{K_{n}( \mu,z_{0},z_{0}) }=e^{i\pi (u-v) }\mathbb{S}(u-v).
\end{gather*}
\end{thm}

In the special case where $\Gamma $ is the unit circle, we have $\Phi (z) =z$, and this reduces to Theo\-rem~\ref{thm4.1}.

\subsection{Bergman polynomials}

Let $G$ be a bounded simply connected domain in the complex plane, bounded by a Jordan curve~$\Gamma $. As above, $\Phi $ is the conformal map of the exterior of $G$ onto the exterior of the unit ball. Let~$\mu $ be a f\/inite positive Borel measure on~$G$. We may def\/ine, for $n\geq 0$, orthonormal polynomials
\begin{gather*}
p_{n}(\mu,z) =\kappa _{n}z^{n}+\cdots,\qquad \kappa _{n}>0
\end{gather*}
satisfying
\begin{gather*}
\int_{G}p_{n}(\mu,z) \overline{p_{m}(\mu,z) }d\mu (z) =\delta _{mn}.
\end{gather*}
We shall assume that $\mu $ is regular in the sense of Stahl and Totik and in this section, also assume that $\mu $ is absolutely continuous with respect to \textit{planar} Lebesgue measure~$dA$ near given points on~$\partial G$. In this sense, the polynomials $\{ p_{n}\} $ fall within the framework of \textit{Bergman polynomials}. For $u\in \partial G$, we def\/ine
\begin{gather*}
\frac{d\mu }{dA}(u) =\lim_{z\rightarrow u,\, z\in G}\frac{d\mu }{dA}(z),
\end{gather*}
whenever the limit is def\/ined. The $n$-th reproducing kernel and its normalized cousin are given by~\eqref{eq4.2} and~\eqref{eq4.3}.

In formulating the result, we need the notion of the convex hull $\operatorname{Co}(K) $ of a set $K$, as well as its boundary $\partial \operatorname{Co}(K)
$. If $J\subset \partial G$, a $\partial G$ neighborhood of $J$ means a~relatively open sub\-set~$J_{1}$ of~$\partial G$ containing~$J$. In \cite{Lubinsky2010}, I proved:

\begin{thm} Let $G$ be a bounded simply connected set, and assume that $\Gamma =\partial G$ is of class $C(1,\alpha) $, with $\alpha \in \left( \frac{1}{2},1\right) $. Let $J\subset \partial G$ be compact, and let some $\partial G$ neighborhood of $J$ also lie in $\partial \operatorname{Co}(G)$.
Let $\mu $ be a finite positive Borel measure on $G$ that is regular. Assume that $\mu $ is absolutely continuous with respect to planar Lebesgue measure in an open subset of $G$ whose boundary contains a $\partial G$ neighborhood of~$J$. Assume moreover, that $\frac{d\mu }{dA}$ is positive and continuous at each point of~$J$. Then uniformly for $z\in J$ and $u$, $v$ in compact subsets of the plane, we have
\begin{gather*}
\lim_{n\rightarrow \infty }\frac{K_{n}\left( \mu,z+\frac{u}{n},z+\frac{v}{n}\right) }{K_{n}( \mu,z,z) }=H\big( u\Phi ^{\prime }(z) \overline{\Phi (z) }+v\overline{\Phi ^{\prime }(z) }\Phi (z) \big),
\end{gather*}
where
\begin{gather*}
H(t) = \begin{cases}
2\dfrac{e^{t}( t-1) +1}{t^{2}}, & t\neq 0, \\
1, & t=0.
\end{cases}
\end{gather*}
\end{thm}

The restriction that $J\subset \partial G\cap \partial \operatorname{Co}(G) $ is a severe geometric restriction~-- basically requiring that~$G$ is ``locally
convex''" in some neighborhood of~$J$. More general boundaries were allowed in~\cite{Lubinsky2010}.

\begin{coro} Let $r$, $s$ be non-negative integers. Then uniformly for $z\in J$,
\begin{gather*}
\lim_{n\rightarrow \infty }\frac{K_{n}^{(r,s) }( \mu,z,z) }{n^{r+s}K_{n}( \mu,z,z) }=\frac{2\big( \Phi
^{\prime }(z) \overline{\Phi (z) }\big) ^{r}\big(\overline{\Phi ^{\prime }(z) }\Phi (z) \big) ^{s}}{r+s+2}.
\end{gather*}
\end{coro}

A substantial generalization of this result has been given by Christopher Sinclair and Maxim Yattselev~\cite{SinclairYattselev2012}. They considered varying, potential theoretic weights, and established Szeg\H{o} asymptotics for the associated orthogonal polynomials, as well as universality limits. There one needs a~generalization of the function~$H$.

\subsection{Rational orthogonal polynomials}

Assume that we are given a sequence of extended complex numbers $\{\alpha _{j}\} \subset \mathbb{\bar{C}}\backslash [-1,1] $, that will serve as our poles. Assume that for some $\eta >0$, and for all $j\geq 1$, $\operatorname{dist}( \alpha _{j},[-1,1]) \geq \eta $. We let $\pi _{0}(x) =1$, and for $k\geq 1$,
\begin{gather*}
\pi _{k}(x) =\prod\limits_{j=1}^{k}( 1-x/\alpha _{j}).
\end{gather*}
Def\/ine nested spaces of rational functions by $\mathcal{L}_{-1}=\{0\}$, $\mathcal{L}_{0}=\mathbb{C}$, and for $k\geq 1,$
\begin{gather*}
\mathcal{L}_{k}=\mathcal{L}_{k}\left\{ \alpha _{1},\alpha _{2},\dots,\alpha
_{k}\right\} =\left\{ \frac{P}{\pi _{k}}\colon \deg (P) \leq k\right\}.
\end{gather*}
Note that if all $\alpha _{j}=\infty $, then $\mathcal{L}_{k}=\mathcal{P}_{k} $. Assume that the poles have an asymptotic distribu\-tion~$\nu$ (with support in $\mathbb{\bar{C}}\backslash [-1,1] $), so that
\begin{gather}\label{eq4.4}
\lim_{k\rightarrow \infty }\log \vert \pi _{k-1}(x)\vert ^{1/k}=\int \log \vert 1-x/t\vert d\nu (t),
\end{gather}
uniformly for $x\in [-1,1] $. Def\/ine orthogonal rational functions $\varphi _{0}$, $\varphi _{1}$, $\varphi _{2}, \dots$ corresponding to the measure $\mu $, such that $\varphi _{k}\in \mathcal{L}_{k}\backslash \mathcal{L}_{k-1}$, and
\begin{gather*}
\int_{-1}^{1}\varphi _{j} \overline{\varphi _{k}} d\mu =\delta _{jk}.
\end{gather*}
Also def\/ine the corresponding rational kernel functions
\begin{gather*}
K_{n}^{r}(\mu,x,y) =\sum_{j=0}^{n-1}\varphi _{j}(x) \overline{\varphi _{j}(y) },
\end{gather*}
and the normalized form
\begin{gather*}
\tilde{K}_{n}^{r}(\mu,x,y) =\mu ^{\prime }(x) ^{1/2}\mu ^{\prime }(y) ^{1/2}K_{n}^{r}( d\mu,x,y).
\end{gather*}
Karl Deckers and the author proved \cite[Theorem~1.2, p.~275]{DeckersLubinsky2012}:

\begin{thm} Let $\mu $ be a regular measure on $[-1,1] $. Let $I$ be an open subinterval of $(-1,1) $ in which $\mu $ is absolutely continuous. Assume that $\mu ^{\prime }$ is positive and continuous at a given \mbox{$x\in I$}. Assume that the poles $\{ \alpha _{j}\} $ are all at least $\eta $ away from $[-1,1] $ and have the asymptotic distribution specified by~\eqref{eq4.4}. Then for $x\in I$ and uniformly for $u$, $v$ in compact subsets of the real line,
\begin{gather*}
\lim_{n\rightarrow \infty }\frac{K_{n}^{r}\left( \mu,x+\frac{u}{\tilde{K}_{n}^{r}(\mu,x,x) },x+\frac{v}{\tilde{K}_{n}^{r}( \mu,x,x) }\right) }{K_{n}^{r}(\mu,x,x) }e^{i\left[ \arg\left( \pi _{n-1}\left( x+\frac{u}{\tilde{K}_{n}^{r}(\mu,x,x) }\right) \right) -\arg \left( \pi _{n-1}\left( x+\frac{v}{\tilde{K}_{n}^{r}(\mu,x,x) }\right) \right) \right] } \\
\qquad{} =\mathbb{S}(u-v).
\end{gather*}
\end{thm}

\subsection{Multivariate orthogonal polynomials}

Let $d\geq 2$, and $\Pi _{n}^{d}$ denote the space of polynomials in $d$ variables of degree at most $n$. Let $N_{n}^{d}=\binom{n+d}{n}$ denote its dimension. (This is a dif\/ferent notion of degree of a multivariate polynomial from that used for Theorem~\ref{thm2.13}.) Let $\mu $ be a positive measure on~$\mathbb{R}^{d}$. We say $\mu $ is \textit{regular}, if
\begin{gather*}
\lim_{n\rightarrow \infty }\left( \sup_{P\in \Pi _{n}^{d}}\frac{\Vert P \Vert _{L_{\infty }( \operatorname{supp}[ \mu]) }^{2}}{\int \vert P\vert ^{2}d\mu }\right) ^{1/n}=1.
\end{gather*}
This is often called the Bernstein--Markov condition \cite{Bloometal2015}, but we prefer the term regularity, to be consistent with the univariate case. We let $K_{n}( \mu,\mathbf{x},\mathbf{y}) $ denote the reproducing kernel for~$\mu $ and~$\Pi _{n}^{d}$, so that\ for all $P\in \Pi _{n}^{d}$, and all $\mathbf{x}\in \mathbb{R}^{d}$,
\begin{gather*}
P( \mathbf{x}) =\int P(\mathbf{y}) K_{n}( \mu,\mathbf{y},\mathbf{x}) d\mu ( \mathbf{y}).
\end{gather*}
Let $J_{\alpha }^{\ast }(z) =z^{-\alpha }J_{\alpha }(z) $ for $\alpha >0$. Andras Kroo and the author proved some general results on universality, which yielded \cite[Theorem~1.7, p.~606]{KrooLubinsky2013}:

\begin{thm} Let $\mu $ be a regular measure on $\bar{B}=\{\mathbf{x}\in \mathbb{R}^{d}\colon \Vert \mathbf{x}\Vert \leq 1\}$, and assume that $D$ is a compact subset of the interior of $\bar{B}$, such that $\mu ^{\prime }$ is positive and continuous in~$D$. Then uniformly for $\mathbf{x}\in D$, and $\mathbf{u}$, $\mathbf{v}$ in compact subsets of~$\mathbb{R}^{d}$,
\begin{gather*}
\lim_{n\rightarrow \infty }\frac{K_{n}\left( \mu,\mathbf{x}+\frac{\mathbf{u}}{n},\mathbf{x}+\frac{\mathbf{v}}{n}\right) }{K_{n}( \mu,\mathbf{x},\mathbf{x}) }=\frac{J_{d/2}^{\ast }\left( \sqrt{G( \mathbf{x},\mathbf{u},\mathbf{v}) }\right) }{J_{d/2}^{\ast }(0) },
\end{gather*}
where if $\cdot $ denotes the standard Euclidean inner product,
\begin{gather*}
G ( \mathbf{x},\mathbf{u},\mathbf{v}) =\Vert \mathbf{u}-\mathbf{v}\Vert ^{2}+\frac{( \mathbf{x}\cdot ( \mathbf{u}-\mathbf{v})) ^{2}}{1-\Vert
\mathbf{x}\Vert ^{2}}.
\end{gather*}
\end{thm}

For the simplex, we have \cite[Theorem~1.8, p.~606]{KrooLubinsky2013}:

\begin{thm} Let $\mu $ be a regular measure on the $d$-dimensional simplex
\begin{gather*}
\Sigma ^{d}=\left\{ \mathbf{x}\in \mathbb{R}^{d}\colon x_{1},x_{2},\dots,x_{d}\geq 0;\sum_{j=1}^{d}x_{j}\leq 1\right\}.
\end{gather*}
 Assume that $D$ is a compact subset of the interior of $\Sigma ^{d}$, such that $\mu ^{\prime }$ is positive and continuous in~$D$. Then uniformly for $\mathbf{x}\in D$, and $\mathbf{u}$, $\mathbf{v}$ in compact subsets of~$\mathbb{R}^{d}$,
\begin{gather*}
\lim_{n\rightarrow \infty }\frac{K_{n}\left( \mu,\mathbf{x}+\frac{\mathbf{u}}{n},\mathbf{x}+\frac{\mathbf{v}}{n}\right) }{K_{n}( \mu,\mathbf{x},\mathbf{x}) }=\frac{J_{d/2}^{\ast }\left( \sqrt{H( \mathbf{x},\mathbf{u},\mathbf{v}) }\right) }{J_{d/2}^{\ast }(0) },
\end{gather*}
where,
\begin{gather*}
H\left( \mathbf{x},\mathbf{u},\mathbf{v}\right) =\sum_{j=1}^{d+1}\frac{(u_{j}-v_{j}) ^{2}}{x_{j}},
\end{gather*}
and the $( d+1) $st component is given by
\begin{gather*}
x_{d+1}=1-\sum_{j=1}^{d}x_{j},
\end{gather*}
with similar definitions of $u_{d+1}$, $v_{d+1}$.
\end{thm}

Universality on the boundary of the $d$-dimensional ball was investigated in~\cite{KrooLubinsky2013B}.

\subsection{Schr\"{o}dinger operators}

Anna Maltsev \cite{Maltsev2010} observed and used analogies between Schr\"{o}dinger operators and orthogonal polynomials to establish universality limits
for Schr\"{o}dinger operators. Let
\begin{gather*}
A=-\frac{d^{2}}{dx^{2}}+V(x)
\end{gather*}
be a Schr\"{o}dinger operator on $L^{2}[0,\infty )$ with Neumann boundary conditions at $x=0$. Assume that $V$ is locally integrable and bounded from below. Let $u( \xi,x) $ be the standard fundamental solution of the eigenvalue equation
\begin{gather*}
Au( \xi,x) =\xi u( \xi,x)
\end{gather*}
with initial conditions
\begin{gather*}
u ( \xi,0 ) =1\qquad \text{and} \qquad u^{\prime } ( \xi,0) =0.
\end{gather*}
If $\mu $ is the spectral measure of $A$, then for $L>0$, its associated reproducing kernel is
\begin{gather*}
S_{L}\left( \xi,\zeta \right) =\int_{0}^{L}u ( \xi,t ) u (\zeta,t) dt,
\end{gather*}
in the sense that
\begin{gather*}
u( \xi,x) \chi _{[0,L] }(x) =\int S_{L}( \xi,\zeta ) u( \zeta,x) d\mu( \zeta).
\end{gather*}
We call a perturbation $q$ \textit{non-destructive} if it leaves the essential spectrum unchanged. We say that it has \textit{zero-average} if
\begin{gather*}
\lim_{L\rightarrow \infty }\frac{1}{L}\int_{0}^{L} \vert q \vert =0.
\end{gather*}
Maltsev proved \cite[p.~464, Theorem~1.3]{Maltsev2010}:

\begin{thm} Let $A=-\frac{d^{2}}{dx^{2}}+p(x) +q(x) $, where $p$ is periodic and continuous, and $q$ is non-destructive and has zero average. Let $d\mu (x)$ be its spectral measure. Let $I$ be a compact subinterval of the essential spectrum, such that~$\mu $ is absolutely continuous on~$I$, while $\mu ^{\prime }$ is continuous and non-zero on~$I$. Let $\xi _{0}\in I$. Then uniformly for $u$, $v$ in bounded subsets of the real line,
\begin{gather*}
\lim_{L\rightarrow \infty }\frac{S_{L}\left( \xi _{0}+\frac{u}{L},\xi _{0}+\frac{v}{L}\right) }{S_{L}\left( \xi _{0},\xi _{0}\right) }=\mathbb{S}(\rho ( \xi _{0}) (u-v)),
\end{gather*}
where $\rho $ is the density of states.
\end{thm}

Maltsev uses this result to study asymptotic spacing of the zeros of~$u$, and also investigated a number of related settings.

\subsection{Universality for entire functions}

Mishko Mitkovski \cite{Mitkovski2013} established a version of universality in the bulk where polynomials are replaced by entire functions. Let $\mu $ be a~positive measure in the real line, which is Poisson summable, that is
\begin{gather*}
\int \frac{d\mu (t) }{1+t^{2}}<\infty.
\end{gather*}
For $T>0$, let $\mathcal{E}_{T} ( \mu ) $ denote the Hilbert space of entire functions of exponential type $\leq T$ that lie in $L_{2} (\mu ) $. It has a~reproducing kernel, which we denote by $K_{T}(\xi,\zeta) $. For complex $z$, def\/ine its majorant
\begin{gather*}
m_{T}(z) =\sup \big\{ \vert F(z) \vert \colon F\in \mathcal{E}_{T}(\mu), \, \Vert F \Vert _{L_{2}(\mu) }\leq 1\big\}.
\end{gather*}
We say that $\mu $ is regular if for any $\varepsilon >0$, there exists $C>0$ such that for all $T>0$ and $t\in \operatorname{supp} [\mu] $,
\begin{gather*}
m_{T}(t) \leq Ce^{\varepsilon T}.
\end{gather*}
Mitkovski \cite[Theorem~1.2]{Mitkovski2013} proved:

\begin{thm} Let $\mu $ be a positive measure on the real line that is Poisson summable, and regular. Assume that $\mu $ is absolutely continuous in a~neighborhood of $\xi _{0}$, and that $\mu ^{\prime}$ is positive and continuous at $\xi _{0}$. Then for all $u,v\in \mathbb{R}$, we have
\begin{gather*}
\lim_{T\rightarrow \infty }\frac{K_{T}\left( \xi _{0}+\frac{u}{T},\xi _{0}+\frac{v}{T}\right) }{K_{T}\left( \xi _{0},\xi _{0}\right) }=\mathbb{S}(u-v).
\end{gather*}
\end{thm}

\subsection{Dirichlet polynomials}

Let $1=\lambda _{1}<\lambda _{2}<\lambda _{3}<\cdots $. Let $\phi _{1}=1$, and for $n\geq 2$, let
\begin{gather*}
\phi _{n}(t) =\frac{\lambda _{n}^{1-it}-\lambda _{n-1}^{1-it}}{\sqrt{\lambda _{n}^{2}-\lambda _{n-1}^{2}}}.
\end{gather*}
One can show that
\begin{gather*}
\int_{-\infty }^{\infty }\phi _{n}(t) \overline{\phi _{m} (t) }\frac{dt}{\pi \left( 1+t^{2}\right) }=\delta _{mn},\qquad m,n\geq 1.
\end{gather*}
The $n$th reproducing kernel for the span of $\{ \lambda_{j}^{-it},\, 1\leq j\leq n\} $ is
\begin{gather*}
K_{n}(x,t) =\sum_{j=1}^{n}\phi _{j}(x) \overline{\phi _{j}(t) }.
\end{gather*}
The author proved \cite[Theorem 1.4]{Lubinsky2014}:

\begin{thm}Assume
\begin{gather*}
\lim_{n\rightarrow \infty }\lambda _{n}=\infty \qquad \text{and} \qquad \lim_{n\rightarrow \infty }\frac{\lambda _{n+1}}{\lambda _{n}}=1.
\end{gather*}
Uniformly for $u$, $v$ in compact subsets of $\mathbb{C}$, and $x$ in compact subsets of the real line,
\begin{gather*}
\lim_{n\rightarrow \infty }\frac{1}{\log \lambda _{n}}K_{n}\left( x+\frac{u}{\log \lambda _{n}},x+\frac{v}{\log \lambda _{n}}\right) =\big[1+x^{2}\big] e^{i(u-v) /2}\mathbb{S}\left( \frac{u-v}{2\pi }\right).
\end{gather*}
\end{thm}

We close this section with three problems. The f\/irst involves generalizations of the results for the unit circle, its subarcs, and smooth
closed curves, discussed in Sections~\ref{section4.1}--\ref{section4.3}:

\begin{problem} Investigate universality limits in the bulk for measures on a finite system of closed curves or arcs in the plane.
\end{problem}

\begin{problem}Investigate universality limits at the edge for measures on a finite system of arcs in the plane.
\end{problem}

Vili Totik noted that this problem is open even for the case of a single smooth arc, other than an interval.

\begin{problem}Investigate universality at interior points of a curve, where a~measure has a~zero or inf\/inity or a discontinuity $($such as a jump$)$. Also, what does universality look like at corners and cusps?
\end{problem}

\section{Varying exponential weights}\label{section5}

The archetypal varying exponential weight is $\exp( -2nx^{2})$, leading to the Gaussian Unitary Ensemble considered by Wigner. The case of general $\exp(-2nQ)$ has been investigated over the decades with varying degrees of rigor. Potential theory plays a crucial role in this endeavor. $Q$ is called an \textit{external field} in this context.

Assume that $\Sigma \subset \mathbb{R}$ is a closed set of positive logarithmic capacity, and if $\Sigma $ is unbounded, that $Q\colon \Sigma \rightarrow \lbrack 0,\infty )$ is continuous, with
\begin{gather*}
\lim_{\vert x\vert \rightarrow \infty,\, x\in \Sigma }\frac{Q(x) }{\log \vert x\vert }=\infty.
\end{gather*}
Associated with $\Sigma $ and $Q$, we may consider the extremal problem
\begin{gather*}
\inf_{\nu }\left( \iint \log \frac{1}{ \vert x-t \vert }d\nu (x) d\nu (t) +2\int Q d\nu \right),
\end{gather*}
where the inf is taken over all positive Borel measures $\nu $ with support in $\Sigma $ and $\nu ( \Sigma) =1$. The inf is attained by a~unique equilibrium measure $\nu _{Q}$, characterized by the following conditions: let
\begin{gather*}
V^{\nu _{Q}}(z) =\int \log \frac{1}{\vert z-t\vert }d\nu _{Q}(t)
\end{gather*}
denote the potential for $\nu _{Q}$. Then \cite[p.~27, Theorem~I.3.1]{SaffTotik1997}
\begin{gather*}
V^{\nu _{Q}}+Q \geq F_{Q}\qquad \text{on}\quad \Sigma ; \\
V^{\nu _{Q}}+Q =F_{Q} \qquad \text{q.e.\ in} \quad \operatorname{supp}[\nu _{Q}].
\end{gather*}
Here the number $F_{Q}$ is a constant, and recall that q.e.\ means except on a set of capacity~$0$.

In the case where $\Sigma $ is an interval, and $Q$ is convex, or $xQ^{\prime }(x) $ exists and is increasing in $( 0,\infty) $ on $\Sigma $, the support of $\nu _{Q}$ is an interval $[ a_{-1},a_{1}] $, called the Mhaskar--Rakhmanov--Saf\/f interval. The Mhaskar--Rakhmanov--Saf\/f numbers $a_{\pm 1}$ are def\/ined by the equations \cite[p.~57, Theorem~2.14]{LevinLubinsky2001}, \cite[p.~201, Theorem~IV.1.11]{SaffTotik1997}
\begin{gather*}
1 =\frac{1}{\pi }\int_{a_{-1}}^{a_{1}}\frac{xQ^{\prime }(x) }{\sqrt{( x-a_{-1}) ( a_{1}-x) }}dx, \qquad
0 =\frac{1}{\pi }\int_{a_{-1}}^{a_{1}}\frac{Q^{\prime }(x) }{\sqrt{( x-a_{-1} )( a_{1}-x) }}dx.
\end{gather*}
The measure $\nu _{Q}$ is absolutely continuous in $(a_{-1},a_{1})$, and its density is given by \cite[p.~42]{LevinLubinsky2001}
\begin{gather}
\nu _{Q}^{\prime }(x) =\frac{\sqrt{( x-a_{-1})(a_{1}-x) }}{\pi ^{2}}\int_{a_{-1}}^{a_{1}}\frac{Q^{\prime }(t) -Q^{\prime }(x) }{t-x}\frac{dt}{\sqrt{(t-a_{-1})( a_{1}-t) }}.\label{eq5.2}
\end{gather}
In the case where $Q$ is even, $a_{-1}=-a_{1}$. See the monographs \cite{Mhaskar1996,SaffTotik1997} or \cite[Chapter~2]{LevinLubinsky2001} for a~comprehensive
introduction.

For the key example
\begin{gather*}
Q(x) =\vert x\vert ^{\alpha },\qquad x\in \mathbb{R},
\end{gather*}
$\alpha >0$, we have $a_{1}=\beta _{a}$, where \cite[p.~204, p.~210]{MhaskarSaff1984}
\begin{gather*}
a_{1}=\left[ \frac{2^{\alpha -2}\Gamma ( \alpha /2) ^{2}}{\Gamma (\alpha) }\right] ^{1/\alpha },
\end{gather*}
and one can determine $\nu _{Q}$ by \eqref{eq5.2}, or by $\nu _{Q}^{\prime }(x) =\nu _{\alpha }( x/a_{1}) /a_{1}$, $x\in (-a_{1},a_{1})$, where \cite[p.~205]{MhaskarSaff1984}
\begin{gather*}
\nu _{\alpha }(x) =\frac{\alpha }{\pi }\int_{\vert t \vert }^{1}\frac{y^{\alpha -1}}{\sqrt{y^{2}-x^{2}}}dx, \qquad x\in (-1,1).
\end{gather*}

The f\/irst rigorous results for general $Q$ were established using the Deift--Zhou steepest descent method, building on the Fokas--Its--Kitaev representation of orthogonal polynomials as solutions of $2\times 2$ matrix Riemann--Hilbert problems. Bleher and Its had earlier~\cite{BleherIts1999} considered the case $Q(x) =x^{4}-tx^{2}$, for a general~$t$. In the mathematical physics literature, the results of Pastur and Shcherbina~\cite{PasturShcherbina1997} are amongst the most general.

Deift, Kriecherbauer, McLaughlin, Venakides, and Zhou \cite{Deiftetal1999,Deiftetal1999B} considered a function~$Q$ (in their terminology $\frac{V}{2}$), that is real analytic on~$\mathbb{R}$, with
\begin{gather}\label{eq5.3}
\lim_{\vert x\vert \rightarrow \infty }\frac{Q(x) }{\log \vert x\vert }=\infty.
\end{gather}
They note that in this case the equilibrium measure $\nu _{Q}$ has support~$J $ consisting of f\/initely many intervals. There is a (mostly) explicit
formula for $\nu _{Q}^{\prime }$, and in particular, $\nu _{Q}^{\prime }$ is positive and analytic in the interior of any of the intervals in the support.

\begin{thm}Let $Q\colon \mathbb{R}\rightarrow \mathbb{R}$ be real valued, and the restriction to the real line of a function analytic in an open set
containing $\mathbb{R}$. Assume that \eqref{eq5.3} holds, and let~$J$ denote the support of the equilibrium measure~$\nu _{Q}$ for the external field $Q$. Then for any $m\geq 1$, $x\in J^{o}$, with $\nu _{Q}^{\prime }(x) >0$, and $u_{1},u_{2},\dots,u_{m}\in \mathbb{R}$, there is the universality limit
\begin{gather*}
\lim_{n\rightarrow \infty }\frac{1}{( n\nu _{Q}^{\prime }(x)) ^{m}}R_{m,n}\left( e^{-2nQ},x+\frac{u_{1}}{n\nu
_{Q}^{\prime }(x) },x+\frac{u_{2}}{n\nu _{Q}^{\prime } (x) },\dots,x+\frac{u_{m}}{n\nu _{Q}^{\prime }(x) }\right) \\
\qquad =\det [ \mathbb{S} ( u_{j}-u_{k} ) ] _{1\leq j,k\leq m}.
\end{gather*}
\end{thm}

This result from \cite[p.~1348, Theorem~1.4]{Deiftetal1999} was a rather direct consequence of far deeper asymptotics for orthogonal polynomials, covering every part of the complex plane. These asymptotics should surely also yield universality at the endpoints of the interval~$J$. However, this was not stated in that paper.

One drawback of the Riemann--Hilbert method is the requirement of analyticity of~$Q$. We shall shortly discuss the $\bar{\partial}$-method, which permits its use for non-analytic~$Q$. However, for universality in the bulk, the most successful general method involved the same mix of methods that yields Theorem~\ref{thm2.5}~-- classical techniques of orthogonal polynomials, complex analysis, and Paley--Wiener space. With the aid of these, and Vili Totik's asymptotics \cite{Totik2000B} for Christof\/fel functions (see also~\cite{Totik2000}), Eli Levin and the author proved~\cite{LevinLubinsky2008}:

\begin{thm}\label{thm5.2} Let $W=e^{-Q}$ be a continuous non-negative function on the set $\Sigma $, which is assumed to consist of at most finitely
many intervals. If $\Sigma $\ is unbounded, we assume also~\eqref{eq5.3}. Let~$h$ be a bounded positive continuous function on~$\Sigma $. Let~$I$ be a closed interval lying in the interior of $\operatorname{supp} [ \nu _{Q}] $, where $\nu _{Q}$ denotes the equilibrium measure for $Q$. Assume that $\nu _{Q}$ is absolutely continuous in a neighborhood of~$I$, and that $\nu_{Q}^{\prime }$ and $Q^{\prime }$ are continuous in that neighborhood, while $\nu _{Q}^{\prime }>0$ there. Then uniformly for $x\in I$, and~$u$,~$v$ in compact subsets of the real line, we have
\begin{gather}\label{eq5.4}
\lim_{n\rightarrow \infty }\frac{\widetilde{K}_{n}\left( hW^{2n},x+\frac{u}{\widetilde{K}_{n}\left( hW^{2n},x,x\right) },x+\frac{v}{\widetilde{K}_{n}\left( hW^{2n},x,x\right) }\right) }{\widetilde{K}_{n}\left(hW^{2n},x,x\right) }=\mathbb{S}(u-v),
\end{gather}
or equivalently,
\begin{gather*}
\lim_{n\rightarrow \infty }\frac{\widetilde{K}_{n}\left( hW^{2n},x+\frac{u}{n\nu _{Q}^{\prime }(x) },x+\frac{v}{n\nu _{Q}^{\prime }(x) }\right) }{n\nu _{Q}^{\prime }(x) }=\mathbb{S}(u-v).
\end{gather*}
\end{thm}

In particular, when $Q^{\prime }$ satisf\/ies a Lipschitz (or if you prefer, H\"{o}lder) condition of some positive order in a neighborhood of~$J$, then \cite[p.~216]{SaffTotik1997} $\nu _{Q}^{\prime }$ is continuous there, and hence we obtain universality except near zeros of $\nu _{Q}^{\prime }$. Theorem~\ref{thm5.2} was a special case of a more implicit result:

\begin{thm}For $n\geq 1$, let $\mu _{n}$ be a positive Borel measure on the real line, with at least the first $2n+1$ power moments finite. Let $I$ be a compact interval in which each $\mu_{n}$ is absolutely continuous. Assume moreover that in~$I$,
\begin{gather*}
d\mu _{n}(x) =h(x) W_{n}^{2n}(x) dx,
\end{gather*}
where
\begin{gather*}
W_{n}=e^{-Q_{n}}
\end{gather*}
is continuous on $I$, and $h$ is a bounded positive continuous function on $I$. Let $\nu _{Q_{n}}$ denote the equilibrium measure for the restriction of $Q_{n}$ to~$I$. Let $J$ be a compact subinterval of $I^{o}$. Assume that
\begin{itemize}\itemsep=0pt
\item[$(a)$]$\{ \nu _{Q_{n}}^{\prime } \} _{n=1}^{\infty }$ are positive and uniformly bounded in some open interval containing~$J$;
\item[$(b)$] $\{ Q_{n}^{\prime }\} _{n=1}^{\infty }$ are equicontinuous and uniformly bounded in some open interval containing~$J$;
\item[$(c)$] for some $C_{1},C_{2}>0$, and for $n\geq 1$ and $x\in I$,
\begin{gather*}
C_{1}\leq K_{n}\big( W_{n}^{2n},x,x\big) W_{n}^{2n}(x) /n\leq C_{2};
\end{gather*}
\item[$(d)$] uniformly for $x\in J$ and $u$ in compact subsets of the real line,
\begin{gather*}
\lim_{n\rightarrow \infty }\frac{K_{n}\left( W_{n}^{2n},x,x\right)
W_{n}^{2n}(x) }{K_{n}\left( W_{n}^{2n},x+\frac{u}{n},x+\frac{u}{n}\right) W_{n}^{2n}\left( x+\frac{u}{n}\right) }=1.
\end{gather*}
Then uniformly for $x\in J$, and $u$, $v$ in compact subsets of the real line, we have~\eqref{eq5.4} with~$hW^{2n}$ replaced by $\mu _{n} $.
\end{itemize}
\end{thm}

Our proof actually established the following limit, uniformly for $x\in J$ and $u$, $v$ in compact subsets of the complex plane, not just the real line:
\begin{gather*}
\lim_{n\rightarrow \infty }\frac{K_{n}\left( \mu _{n},x+\frac{u}{\widetilde{K}_{n}\left( \mu _{n},x,x\right) },x+\frac{v}{\widetilde{K}_{n}(\mu
_{n},x,x) }\right) }{K_{n} ( \mu _{n},x,x ) }e^{-\frac{n}{\widetilde{K}_{n} ( \mu _{n},x,x ) }Q_{n}^{\prime }(x)( u+v) }=\mathbb{S}(u-v).
\end{gather*}

In a recent paper \cite{LevinLubinsky2015}, Eli Levin and the author established that universality holds in measure for sequences of varying weights under hypotheses weaker than those in Theorem~\ref{thm5.2}. In principle, the same sorts of techniques that yield Theorem~\ref{thm5.2} in the bulk,
should work at the soft and hard edges. Indeed, the equivalence of universality along the diagonal (that is the limit~\eqref{eq5.4} with $a=b$) to universality for $u$, $v$ in compact sets, was explored for the soft edge in~\cite{LevinLubinsky2010} and the hard edge in~\cite{Lubinsky2008D}. However, the asymptotics for Christof\/fel functions that are needed to apply these results are not generally available.

The most general known universality results at the edge are due to McLaughlin and Miller~\cite{McLaughlinMiller2008}. They use the Riemann--Hilbert method, but modif\/ied using $\bar{\partial}$-techniques to approximate non-analytic~$Q$ with analytic~$Q$. They consider the case of an equilibrium measure supported on f\/initely many intervals, and place conditions on the sign of certain functions formed from complex potentials and equilibrium densities:

\begin{thm} Let $Q\colon \mathbb{R}\rightarrow \mathbb{R}$ and assume that $Q^{\prime \prime }$ satisfies a Lipschitz condition of order~$1$. Assume that the support of the equilibrium measure $\nu _{Q}$ consists of finitely many intervals $[ \alpha _{j},\beta _{j}] $, $1\leq j\leq \ell $. Assume that Condition~$2$ in {\rm\cite[p.~17]{McLaughlinMiller2008}} is satisfied. Then at every edge $\xi \in \{ \alpha _{1},\beta _{1},\dots,\alpha _{\ell },\beta _{\ell }\} $, there is a number $\lambda $ such that for
real $u$, $v$,
\begin{gather*}
\lim_{n\rightarrow \infty }\frac{1}{(\lambda n) ^{2/3}}\tilde{K}_{n}\left( e^{-2nQ},\xi +\frac{u}{(\lambda n) ^{2/3}},\xi +\frac{v}{(\lambda n) ^{2/3}}\right) =\mathbb{A}{\rm i}(u,v).
\end{gather*}
In particular, Condition~$2$ is satisfied if~$Q$ is strictly convex, and grows faster at $\infty $ than $( \log\vert x\vert ) ^{1+\varepsilon }$, for some $\varepsilon >0$.
\end{thm}

McLaughlin and Miller also established universality in the bulk of course~-- and both universality limits are consequences of far deeper asymptotics for orthogonal polynomials that hold everywhere in the complex plane. For varying exponential weights on the unit circle~\cite{McLaughlinMiller2006}, McLaughlin and Miller also established asymptotics of the associated orthogonal polynomials that undoubtedly imply universality limits in the bulk, and probably at appropriate edges too.

The Deift--Zhou steepest descent method has also been used to great ef\/fect by Baik, Kriecher\-bauer, McLaughlin, and Miller~\cite{Baiketal2007} in establishing universality associated with sequences of discrete measures $\{ \mu _{n}\} $, where each $\mu _{n}$ has f\/initely many jumps. Appropriate assumptions are placed on the distributions of the jumps, and their size. Universality at the edge or in the bulk (appropriately interpreted), is a consequence of far deeper results on asymptotics of orthogonal polynomials.

Another important recent development is ``global asymptotics'' due to Kriecherbauer, Schubert, Sch\"{u}ler, and Venker \cite{Kriecherbaueretal2015}, where they obtain universality with error estimates, that are uniform in the range, as well as (remarkably) uniform for $Q$ in a class of real analytic external f\/ields. Both universality in the bulk and at the edge are established. See also~\cite{KriecherbauerVenker2016}.

For varying exponential weights in the plane, of the form $e^{-n(\vert z\vert ^{2}-\operatorname{Re}(tz^{2}))}$ on the complex plane, Roman Riser investigated universality in~\cite{Riser2013}, using orthogonal polynomial techniques. There the limiting kernel is not a sine kernel. In a related vein, Antti Haimi~\cite{Haimi2014} established universality for weights $e^{-nQ}$ def\/ined on the whole complex plane, and polyanalytic polynomials, leading to universality limits that involve associated Laguerre polynomials. Polyanalytic Ginibre ensembles were considered, for example, in~\cite{HaimiHedenmalm2013}.

\section{Fixed exponential weights}\label{section6}

A powerful 1999 paper of Deift, Kriecherbauer, McLaughlin, Venakides, and Zhou \cite{Deiftetal1999B} establishes asymptotics in all regions of the
plane for orthogonal polynomials associated with f\/ixed exponential weights $e^{-2Q}$ on $\mathbb{R}$, where $Q$ is a polynomial of even degree with positive leading coef\/f\/icient. Using direct substitution into the Christof\/fel--Darboux formula, this must lead to universality both in the bulk
and at the soft-edge, although this was not stated there. Another powerful paper that implies universality results, but which are not explicitly stated, is that of Kriecherbauer and McLaughlin \cite{KriecherbauerMcLaughlin1999}, for $\exp( -\vert x\vert^{\alpha })$, all $\alpha >0$. Of course, there was a lot of earlier work for such asymptotics, that also implies universality in the bulk, at least, for various exponential weights.

An excellent illustration of universality in the bulk, and at the soft and hard edge, is provided by the 2007 paper of Vanlessen~\cite{Vanlessen2007} for generalized Laguerre weights,
\begin{gather}\label{eq6.1}
\mu ^{\prime }(x) =x^{\alpha }e^{-Q(x) },\qquad x\in (0,\infty ),
\end{gather}
where $\alpha >-1$ and
\begin{gather}\label{eq6.2}
Q(x) =\sum_{k=0}^{m}q_{k}x^{k},
\qquad \mbox{where} \quad m\geq 1 \quad \text{and} \quad q_{m}>0.
\end{gather}
Because we are dealing with a f\/ixed weight, there is a sequence of equilibrium densities and a sequence of Mhaskar--Rakhmanov--Saf\/f intervals $[0,\beta _{n}] $. The $n$th Mhaskar--Rakhmanov--Saf\/f number $\beta_{n}$ is def\/ined by the equation
\begin{gather*}
\frac{1}{2\pi }\int_{0}^{\beta _{n}}Q^{\prime }(x) \sqrt{\frac{x}{\beta _{n}-x}}dx=n.
\end{gather*}
$\beta _{n}$ grows like $n^{1/m}$ and has a complete asymptotic expansion
\begin{gather*}
\beta _{n}=n^{1/m}\left( \beta ^{(0) }+\sum_{k=1}^{\infty }\beta^{( k) }n^{-k/m}\right).
\end{gather*}
The equilibrium density for the rescaled external f\/ield $\frac{1}{n}Q(\beta _{n}x) $ has the form
\begin{gather*}
\nu _{n}^{\prime }(x) =\frac{1}{2\pi }\sqrt{\frac{1-x}{x}}h_{n}(x), \qquad x\in [ 0,1 ],
\end{gather*}
where $h_{n}$ is a polynomial of degree $m-1$, that converges to a~limit polynomial~$h$ as $n\rightarrow \infty $. We can now state the associated
universality limits in contracted form, for the bulk, soft edge at $\beta_{n}$, and hard edge at~$0$~\cite{Vanlessen2007}:

\begin{thm}Assume that $\mu $ is given by \eqref{eq6.1} and \eqref{eq6.2}, and that $\{ \beta _{n}\}$, $\{ \nu _{n}^{\prime }\}$ are as above.
\begin{itemize}\itemsep=0pt
\item[$(a)$] Uniformly for $x$ in compact subsets of $(0,1) $ and $u$, $v$ in compact subsets of the real line,
\begin{gather*}
\frac{\beta _{n}}{n\nu _{n}^{\prime }(x) }\tilde{K}_{n}\left(\mu,\beta _{n}\left( x+\frac{u}{n\nu _{n}^{\prime }(x) }\right)
,\beta _{n}\left( x+\frac{v}{n\nu _{n}^{\prime }(x) }\right) \right) =\mathbb{S}(u-v) +O\left( \frac{1}{n}\right).
\end{gather*}
\item[$(b)$] For $n\geq 1$, let
\begin{gather*}
c_{n}=\left( \frac{1}{2}h_{n}(1) \right) ^{2/3}.
\end{gather*}
 Uniformly for $u$, $v$ in compact subsets of $\mathbb{R}$,
\begin{gather*}
\frac{\beta _{n}}{c_{n}n^{2/3}}\tilde{K}_{n}\left( \mu,\beta _{n}\left( x+\frac{u}{c_{n}n^{2/3}}\right),\beta _{n}\left(x+\frac{v}{c_{n}n^{2/3}}
\right) \right) =\mathbb{A}{\rm i}(u,v) +O\left( \frac{1}{n^{1/3}}\right).
\end{gather*}
\item[$(c)$] For $n\geq 1$, let
\begin{gather*}
\tilde{c}_{n}=\left( \frac{1}{2}h_{n}(0) \right) ^{2}.
\end{gather*}
Uniformly for $u$, $v$ in bounded subsets of $(0,\infty) $,
\begin{gather*}
\frac{\beta _{n}}{4\tilde{c}_{n}n^{2}}\tilde{K}_{n}\left( \mu,\frac{u}{4\tilde{c}_{n}n^{2}},\frac{v}{4\tilde{c}_{n}n^{2}}\right) =\mathbb{J}_{\alpha }(u,v) +O\left( \frac{u^{\alpha /2}v^{\alpha /2}}{n}\right).
\end{gather*}
\end{itemize}
\end{thm}

\looseness=-1 In \cite{LevinLubinsky2009}, Eli Levin and the author used f\/irst order asymptotics for orthogonal polyno\-mials~$\{ p_{n}\} $ for exponential weights, established in~\cite{LevinLubinsky2001}, to obtain universality limits for a broad class of exponential weights $e^{-2Q}$ on the whole real line. Here $Q$ is assumed twice dif\/ferentiable, and satisf\/ies some other regularity conditions, and we also considered weights~$he^{-2Q}$, where $h$ does not oscillate or grow too rapidly~-- for example it could be a~generalized Jacobi weight. The main observation, that all one needs is asymptotics with only a $o(1) $ error term, was due to Eli Levin. In our later paper \cite[p.~720~f\/f.]{LevinLubinsky2008} (ironically published earlier), Eli Levin and the author turned results for varying exponential weights into ones for a more general class of exponential weights:

\begin{dfn} Let $I= ( c,d ) $ be an open interval, bounded or unbounded, containing $0$ in its interior. Let $W=\exp (-Q) $, where $Q\colon I\rightarrow \lbrack 0,\infty )$ satisf\/ies the following properties:
\begin{itemize}\itemsep=0pt
\item[(a)] $Q^{\prime }$ is continuous in $I$ and $Q(0) =0$.
\item[(b)]$Q^{\prime }$ is non-decreasing in $I$;
\item[(c)]
\begin{gather*}
\lim_{t\rightarrow c+}Q(t) =\lim_{t\rightarrow d-}Q(t) =\infty.
\end{gather*}
\item[(d)] The function
\begin{gather*}
T(t) =\frac{tQ^{\prime }(t) }{Q(t) },\qquad t\neq 0,
\end{gather*}
is quasi-increasing in $( 0,d) $, in the sense that for some $C>0$,
\begin{gather*}
0<x<y\Rightarrow T(x) \leq CT(y).
\end{gather*}
$T$ is also assumed quasi-decreasing in $(c,0) $. In addition, we assume that for some $\Lambda >1$,
\begin{gather*}
T(t) \geq \Lambda \qquad \text{in} \quad I\backslash \{0\}.
\end{gather*}
\item[(e)] There exists $\varepsilon _{0}\in (0,1)$, $C_{1},C_{2}>0$ such that for $y\in I\backslash \{0\} $,
\begin{gather*}
C_{1}\leq T(y) /T\left( y\left[ 1-\frac{\varepsilon _{0}}{T(y) }\right] \right) \leq C_{2}.
\end{gather*}
\item[(f)] For every $\varepsilon >0$, there exists $\delta >0$ such that for all $x\in I\backslash \{0\} $,
\begin{gather*}
\int_{x-\frac{\delta \vert x\vert }{T(x) }}^{x+\frac{\delta \vert x\vert }{T(x) }}\frac{Q^{\prime } (s) -Q^{\prime }(x) }{s-x}ds\leq \varepsilon \vert Q^{\prime }(x) \vert.
\end{gather*}
Then we write $W\in \mathcal{F} ({\rm dini} ) $.
\end{itemize}
\end{dfn}

These conditions, especially that in (f) are somewhat technical! One explicit condition on $Q$ that guarantees that those in (e), (f) hold, is
\begin{itemize}\itemsep=0pt
\item[
(g)] $Q^{\prime \prime }$ exists in $I\backslash \{ 0\} $ and there exists $C_{1}>0$ such that
\begin{gather*}
\frac{Q^{\prime \prime }(x) }{\vert Q^{\prime } (x) \vert }\leq C_{1}\frac{Q^{\prime }(x) }{Q(x) } \qquad \text{a.e.} \quad x\in \mathbb{R}\backslash \{0\}.
\end{gather*}
\end{itemize}

Examples of weights in this class are $W=\exp ( -Q ) $ on $I=\mathbb{R}$, where
\begin{gather*}
Q(x) = \begin{cases}
Ax^{\alpha }, & x\in \lbrack 0,\infty ), \\
B\vert x\vert ^{\beta }, & x\in ( -\infty,0),
\end{cases}
\end{gather*}
with $\alpha,\beta >1$. More generally, if $\exp _{k}=\exp ( \exp ( \cdots \exp ( {\,} ) ) ) $ denotes the $k$th iterated exponential, we may take
\begin{gather*}
Q(x) =
\begin{cases}
\exp _{k}\left( Ax^{\alpha }\right) -\exp _{k}(0),
& x\in \lbrack 0,\infty ), \\
\exp _{\ell }\left( B\vert x\vert ^{\beta }\right) -\exp _{\ell
}(0), & x\in ( -\infty,0 ),
\end{cases}
\end{gather*}
where $k,\ell \geq 0$, $\alpha,\beta >1$.

A key descriptive role for such $Q$ is played by the Mhaskar--Rakhmanov--Saf\/f numbers $a_{-n}<0<a_{n}$, def\/ined for $n\geq 1$ by the equations
\begin{gather}
n =\frac{1}{\pi }\int_{a_{-n}}^{a_{n}}\frac{xQ^{\prime }(x) }{\sqrt{( x-a_{-n}) ( a_{n}-x) }}dx,\label{eq6.3} \\
0 =\frac{1}{\pi }\int_{a_{-n}}^{a_{n}}\frac{Q^{\prime }(x) }{\sqrt{( x-a_{-n} ) ( a_{n}-x) }}dx.\label{eq6.4}
\end{gather}
In the case where $Q$ is even, $a_{-n}=-a_{n}$. Yes, there is a conf\/lict with the earlier notation for recurrence coef\/f\/icients~-- in this section $a_{n}$ has a dif\/ferent meaning. We also def\/ine,
\begin{gather}\label{eq6.5}
\delta _{n}=\frac{1}{2} ( a_{n}+ \vert a_{-n} \vert ).
\end{gather}
We proved \cite[Theorem~7.4, p.~771]{LevinLubinsky2008}:

\begin{thm} Let $W=\exp (-Q) \in \mathcal{F}({\rm dini})$. Let $0<\varepsilon <1$. Then uniformly for $u$, $v$ in compact subsets of the real line, and $x\in [ a_{-n}+\varepsilon \delta _{n},a_{n}-\varepsilon \delta _{n} ] $, we have
\begin{gather*}
\lim_{n\rightarrow \infty }\frac{\tilde{K}_{n}\left( W^{2},x+\frac{u}{\tilde{K}_{n}\left( W^{2},x,x\right) },x+\frac{v}{\tilde{K}_{n}\left(
W^{2},x,x\right) }\right) }{\tilde{K}_{n}\left( W^{2},x,x\right) }=\mathbb{S}(u-v).
\end{gather*}
In particular, if $W$ is even, this holds uniformly for $\vert x\vert \leq ( 1-\varepsilon) a_{n}$.
\end{thm}

As far as the author is aware, there are no results on universality at the edge anywhere as general as those in the bulk. Perhaps still the most general result is due to Deift and Gioev~\cite{DeiftGioev2007}. Undoubtedly the varying weights results of McLaughlin and Miller~\cite{McLaughlinMiller2008} imply universality at the edge for appropriate f\/ixed exponential weights, but this does not seem to have been written down. So we close with the result of Deift and Gioev, which was part of a more general treatment of edge universality for orthogonal, symplectic, and unitary ensembles:

\begin{thm}\label{thm6.4} Let $Q$ be a polynomial of positive even degree, with positive leading coefficient. Let $W=e^{-2Q}$ on~$\mathbb{R}$.
Let $a_{\pm n}$ be the $n$th Mhaskar--Rakhmanov--Saff number, defined by~\eqref{eq6.3},~\eqref{eq6.4}. Let $\delta _{n}$ be defined by~\eqref{eq6.5}. Then uniformly for~$u$,~$v$ in compact subsets of $\mathbb{R}$,
\begin{gather*}
\lim_{n\rightarrow \infty }\frac{\delta _{n}}{\tau _{n}n^{2/3}}\tilde{K}_{n}\left( \mu,a_{n}+\frac{\delta _{n}u}{\tau _{n}n^{2/3}},a_{n}+\frac{\delta _{n}v}{\tau _{n}n^{2/3}}\right) =\mathbb{A}{\rm i}(u,v).
\end{gather*}
Here $\{ \tau _{n}\} $ is a sequence of numbers depending on $Q$ and arising from certain equilibrium densities.
\end{thm}

\begin{problem} Establish universality at the edge for more general fixed exponential weights than those in Theorem~{\rm \ref{thm6.4}}.
\end{problem}

\subsection*{Acknowledgements}

Research supported by NSF grant DMS1362208.

It was a privilege to attend the very high level conference celebrating Percy Deift's 70th birthday. I owe Percy a great deal: it was Percy's 60th birthday conference at the Courant Institute that inspired me to try apply classical methods of orthogonal polynomials to universality limits. Percy's comments and perspectives, have really helped in this endeavor. Thank you, Percy~-- and thank you to CRM and the organizers of the conference.

This survey has benef\/ited greatly from the corrections and comments of Thomas Bothner, Jonathan Breuer, Tivadar Danka, Thomas Kriecherbauer, Arno Kuijlaars, Anna Maltsev, Andrei Mart{\'{\i}}nez-Finkelshtein, Barry Simon, Vili Totik, Yu-Qiu Zhao, and the anonymous referees.


\pdfbookmark[1]{References}{ref}
\LastPageEnding

\end{document}